\documentclass[a4paper,12pt,oneside,onecolumn,number]{article}

\usepackage{amsfonts}
\usepackage{amsmath}
\usepackage{amssymb}
\usepackage{amsthm}
\usepackage{cite}
\usepackage{enumitem}
\usepackage[top=20mm,bottom=30mm,left=20mm,right=20mm]{geometry}

\numberwithin{equation}{section}
\newtheorem{theorem}{Theorem}[section]
\newtheorem{corollary}{Corollary}[section]
\newtheorem{remark}{Remark}[section]
\newtheorem{lemma}{Lemma}[section]

\newtheorem{prop}[theorem]{Proposition}

\newcommand{\etal}{\emph{et al}.}

\newcommand{\R}{\mathbb{R}}

\newcommand{\N}{\mathbb{N}}
\newcommand{\C}{\mathbb{C}}

\title{On Uniqueness Theorems for the Inverse $q$-Sturm–Liouville problems}
\author{
F.\,A.~Gawish $^{1}$\thanks{e-mail: fatma.gawish@fsc.bu.edu.eg},
Z.\,S.~Mansour$^{2}$\thanks{e-mail: zsmansour@cu.edu.eg},
}

\begin{document}
\date{\empty}
\maketitle
\begin{center}
\small{$^{1}$\footnotesize{\it Department of Mathematics, Faculty of Science, Benha University, Benha-Kalubia 13518, Egypt}}
\small{$^{2}$\footnotesize{\it Department of Mathematics, Faculty of Science, Cairo University, Giza 12613, Egypt}}\\
\end{center}

\begin{abstract}
We establish two theorems that illustrate the uniqueness of inverse $q$-Sturm-Liouville problems based on a specified set of spectral data. The first uniqueness theorem employs the method of transformation operators to provide a $q$-analog of the Levinson-Marchenko uniqueness theorem.
The second uniqueness theorem is a $q$-analog of the Ashrafyan uniqueness theorem. We introduced a $q$-analog of the Gelfand-Levitan approach, which involves converting $q$-difference operators into $q$-integral operators to prove the second uniqueness theorem.\smallskip\\

\end{abstract}

\section{Introduction}
\par Let $[a, b] \subseteq \mathbb{R}$ be a finite closed interval and $v(\cdot)$ be a continuous real-valued function defined on $[a, b]$. The Sturm-Liouville problem is the problem of finding a function $y(\cdot)$ and a number $\lambda \in \mathbb{C}$ satisfying the differential equation
\begin{equation}\label{sl1}
  l y:=-y^{\prime \prime}+v(x) y(x)=\lambda y(x), \quad a \leqslant x \leqslant b,
\end{equation}
together with the boundary conditions
\begin{align}\label{sl2}
a_{1} y(a)+a_{2} y^{\prime}(a)=0,\quad 
 b_{1} y(b)+b_{2} y^{\prime}(b)=0.
\end{align}
The sequence of eigenvalues $\left\{\lambda_{n}\right\}_{n=0}^{\infty}$ is called the spectrum of $l$ corresponding to the eigenfunctions $\left\{\theta(\cdot,\lambda_n)\right\}_{n=0}^{\infty}$. 
 Recovering operators from their spectral values is the inverse problem of spectral analysis. In numerous domains of mathematics and its applications, such as physics, engineering, and control theory, inverse problems are crucial. For instance, in \cite{casti80}, Casti introduced a general inverse problem in optimal control theory. In~\cite{FERLAUTO06}, Ferlauto and Marsilio discussed numerical techniques to solve two-dimensional inverse problems related to aerodynamic design. Furthermore, Alhaidari and Taiwo in \cite{alhaidari22} focused on determining the potential function linked to the energy spectrum.
 
 In the following, we use the notation $l(v,h,H)$ to donate the operator defined in \eqref{sl1} with the boundary conditions
\begin{equation}\label{slp1}
y^{\prime}(0)-h y(0)=0, 
  \end{equation}
  \begin{equation}\label{slp2}
 y^{\prime}(\pi)+H y(\pi)=0,
  \end{equation}
where $h$ and $H$ are real numbers. Moreover, the norming constants $ \{ a_{n} \}_{n=0}^\infty$ are  defined by 
\begin{equation*}\label{is33}
      a_{n}=\int_{0}^{\pi}|\theta(x,\lambda_{n})|^2 dx,\quad n \in \N_{0},
 \end{equation*}
 where $\{\theta(x,\lambda_{n})\}_{n=0}^\infty$ is the set of eigenfunctions of the problem~\eqref{sl1}-\eqref{slp1}.
 The operator $\tilde{l}$ is defined in a similar manner to the operator $l$, but with the replacements of $v$, $h$, and $H$ by $\tilde{v}$, $\tilde{h}$, and $\tilde{H}$, respectively.
 
In 1949, Ambarzumian introduced the following uniqueness theorem, which is considered the first uniqueness theorem in inverse problems; see~\cite{ambarzumian29}.

\vskip .5 cm

\noindent{\bf Theorem A. }
(\noindent{\bf{Ambarzumian's Theorem}})  Let $\{\lambda_{n}\}_{n=0}^\infty$ denote the eigenvalues of the Sturm-Liouville problem
 \begin{equation*}\label{Am1}
   y^{\prime \prime}+\{\lambda-v(x)\} y =0 \quad(0 \leqslant x \leqslant \pi),
 \end{equation*}
with the boundary conditions
 \begin{equation*}\label{Am2}
   y^{\prime}(0)=y^{\prime}(\pi)=0,
 \end{equation*}
where $v(x)$ is a real-valued continuous function on the interval $[0, \pi]$. If $\lambda_{n}=n^{2}, n\in \N_{0}$, then $v(x)= 0$ a.e. on $(0, \pi).$  

\vskip .5 cm 

Ambarzumian's Theorem is an oddity, and Borg in \cite{borg46} proved that the specification of two spectra of Sturm-Liouville operators uniquely defines the potential function $v(\cdot)$.  In 1949, Levinson presented a different proof of Borg's result.

Levinson was the first to use the contour integral method to study the inverse problem for the Sturm-Liouville operator~\cite{levinson49}. He introduced Theorem B below. Tikhonov (1950) proved the uniqueness theorem for solving the inverse Sturm-Liouville problem on the half-line by utilizing the Weyl function \cite{tikhonov49}. The transformation operator played a significant role in the spectral theory of Sturm-Liouville operators. Marchenko initially employed the transformation operators to introduce another proof for Theorem B, see also~\cite{marchenko52}.
\vskip .5 cm 

\noindent{\bf Theorem B.} (\textbf{Levinson-Marchenko Theorem})

Let $l(v,h,H)$ and $\tilde{l}=l(\tilde{v},\tilde{h},\tilde{H})$ be the operators associated with the spectral data $\{\lambda_{n}, a_{n}\}$ and $\{\tilde{\lambda_{n}},\tilde{a_{n}}\},$ respectively. If $\lambda_{n}=\tilde{\lambda}_{n}, a_{n}=\tilde{a}_{n}, n \in \N_{0}$, then $l=\tilde{l}$, i.e., $v(x)=\tilde{v}(x)$ a.e. on $ (0, \pi), h=\tilde{h} $ and $ H=\tilde{H} $.

\vskip .5 cm 
The isospectrality problem involves characterizing all problems of the form \eqref{sl1}-\eqref{sl2} that have the same spectrum. I.e., two problems are isospectral if they have the same set of eigenvalues. This problem has been studied in detail by Trubowitz et al. in three papers \cite{Turb1,Turb2,Turb3} and by Jodiet and Levitan in \cite{jodeit97}. In \cite{ashrafyan17}, Ashrafyan proved Theorem C below, which is a generalization of the Levinson-Marchenko uniqueness theorem.

\vskip .5 cm

\noindent{\bf Theorem C. }  $(\nonscript{\bf{Ashrafyan\ Theorem.}})$ 

  Let $l(v,h,H)$ and $\tilde{l}=l(\tilde{v},\tilde{h},\tilde{H})$ be the operators associated with the spectral data $\{ \lambda_{n}, a_{n} \}$ and $\{ \tilde{\lambda_{n}},\tilde{a_{n}} \},$ respectively. If $\lambda_{n}=\tilde{\lambda}_{n}, a_{n}\geq \tilde{a}_{n}, n\in \N_{0}$, then $l=\tilde{l}$, i.e. $v(x)=\tilde{v}(x)$ a.e. on $(0, \pi), h=\tilde{h}$ and $H=\tilde{H}$.

\vskip .5 cm 

\par This paper has the following structure after this introduction section, Section 1. Section~2 contains the q-notations that are essential to our task. In Section~3, we introduce the $q$-analog of the Levinson-Marchenko uniqueness theorem, Theorem B. Section~4 includes a $q$-analog of the Gelfand-Levitan equation with some essential properties of $q$-linear operators that we need in Section~5. In Section~5, we introduce the $q$-analog of the Ashrafyan uniqueness theorem, Theorem C. 

\section{$q$-Notation and results}

 This section outlines the fundamental $q$-notations employed in our study. Unless specified otherwise, $q$ is a positive number within the interval $(0, 1)$, $\N$ denotes the set of positive integers, and $\N_0 = \N \cup \{0\}$.
The $q$-shifted factorial is defined for $a \in \mathbb{C}$ by
\begin{equation*}
(a ; q)_n=\left\{\begin{array}{cc}
1, & n=0, \\
\prod_{i=0}^{n-1}\left(1-a q^i\right), & n \in \mathbb{N}.
\end{array}\right.
\end{equation*}
  
  In \cite{jackson03},  Jackson defined a $q$-analog of the exponential function $ E_{q}(z)$ by
  
\begin{equation*}
  E_q(z):=\sum_{n=0}^{\infty}  \frac{q^{n(n-1)/2}(z(1-q))^n}{(q;q)_n}, \quad z \in \mathbb{C}.
\end{equation*}
 A pair of $q$-analogs of the trigonometric functions $\cos(z)$ and $\sin(z)$ are defined on $\C$ by
  \begin{equation}\label{cos}
  \begin{split}
    \cos (z ; q) & :=\sum_{n=0}^{\infty}(-1)^n q^{n^2} \frac{(z(1-q))^{2 n}}{(q ; q)_{2 n}},\\
\sin(z ; q) &:= \sum_{n=0}^{\infty}(-1)^n q^{n^2+n} \frac{(z(1-q))^{2 n+1}}{(q ; q)_{2 n+1}}.
 \end{split}
 \end{equation}
  The trigonometric functions $\cos z$ and $\sin z$ have several $q$-analogs in the literature. See, for example~\cite{gasper2004,koekoek96,askey99}. In \cite{koelink94}, Koelink and Swarttouw demonstrated that the zeros of $\cos(z; q)$ and $\sin(z; q)$ are real and simple.  In \cite{annaby11}, the authors introduced an asymptotic formula for the zeros of $\cos(z; q)$ and $\sin(z; q)$.  

 A $q$-geometric set $A$ is the set such that $q x \in A$ whenever $x \in A$, see~\cite{p2}. Let $f$ be a function, real or complex-valued, defined on a $q$-geometric set $A$, where $|q| \neq 1$.

Jackson in \cite{jackson09} defined the $q$-difference operator, $D_q$, by 
\begin{equation*}\label{}
  D_{q} f(x):=\frac{f(x)-f(q x)}{x-q x},\quad x \in A \backslash\{0\}.
\end{equation*}
 If $0 \in A$, the $q$-derivative at zero is defined for $|q|<1$ by
 \begin{equation*}\label{}
  D_{q} f(0):=\lim _{n \rightarrow \infty} \frac{f\left(x q^{n}\right)-f(0)}{x q^{n}},\quad x \in A \backslash\{0\},
 \end{equation*}
provided the limit exists and is independent of $x$. Furthermore, the $q$-derivative at zero is defined for $|q|>1$ by
$$
D_{q^{-1}} f(0):=D_{q} f(0).
$$
The $q$-product rule is 
\begin{equation*}\label{}
  D_{q}(f g)(x)=g(x) D_{q} f(x)+f(q x) D_{q} g(x),
\end{equation*}
where $f,\ g$ are defined on a $q$-geometric set $A$ and $x\in A$. 
 In \cite{jackson10}, Jackson introduced a $q$-extension of the Riemann integral by
\begin{equation}\label{int}
  \int_{0}^{x} f(t)\,  d_{q} t:=(1-q) \sum_{n=0}^{\infty} x q^{n} f\left(x q^{n}\right) \quad(x \in A),
\end{equation}
provided that the series converges. In general,
$$
\int_{a}^{b} f(t)\, d_{q} t:=\int_{0}^{b} f(t)\, d_{q} t-\int_{0}^{a} f(t)\, d_{q} t, \quad a, b \in A.
$$
   A function $f$ defined on a $q$-geometric set $A$ is said to be $q$-regular at zero if
$$ \lim_{n\rightarrow\infty} f\left(x q^n\right)=f(0),\ \text{ for all } \ x \in A,$$ see~\cite{annabyq12}.
 The fundamental theorem of $q$-calculus states that if $f$ is a function defined on a $q$-geometric set $A$, then for $x\in A$
   \begin{equation}\label{ibp3}
     \quad D_{q} \int_{0}^{x} f(t)\, d_{q} t=f(x), \quad \int_{0}^{x} D_{q} f(t)\, d_{q} t=f(x)-\lim _{n \rightarrow \infty} f\left(x q^{n}\right).
   \end{equation}
   Moreover, if $f$ is $q$-regular at zero, then for $\{a,b\}\subseteq A,$
   \begin{equation*}\label{}
     \quad  \int_{a}^{b}D_{q} f(t)\, d_{q} t=f(b)-f(a),
   \end{equation*}
   see~\cite{annabyq12,gasper2004}. Similarly, in \cite{annabyqq12}, the authors proved that if $f$ is a function defined on a q-geometric set A and $x \in q A$, then
\begin{align}\label{ibp4}
\nonumber & D_{q^{-1}} \int_0^x f(t)\, d_q t=q f\left(x/q\right),\\
& \int_0^x D_{q^{-1}} f(t)\, d_q t=q \left(f(x/q)-\lim _{n \rightarrow \infty} f(x q^{n-1})\right).
\end{align}
Moreover, if $f$ is $q$-regular at zero, then for $\{a,b\}\subseteq A,$

\begin{equation*}\label{}
  \int_a^b  D_{q^{-1}} f(t)\, d_q t=q \left(f(b/q)\right)-f\left(a/q)\right).
\end{equation*}

  Let $f(x,t)$ be a function defined on $A\times A$, then simple manipulation yields

 \begin{equation}\label{ff1}
   D_{q^{-1},x} \int_{0}^{x} f(x,t) \,d_{q} t= \int_{0}^{x} D_{q^{-1},x} f(x,t)\, d_{q} t+  f\left(\frac{x}{q},\frac{x}{q}\right),
 \end{equation}
 \begin{equation}\label{ff2}
   D_{q,x}\int_{0}^{x}  D_{q^{-1},x}f(x,t) \,d_{q} t= \int_{0}^{x}  D_{q,x}D_{q^{-1},x} f(x,t) \,d_{q} t+  D_{q,x}f\left(x,t\right)\mid_{t=x}.
 \end{equation}
The   rules  \eqref{ibp1} and \eqref{ibp2}  below of $q$-integration by parts are derived from the fundamental theorems of $q$-calculus \eqref{ibp3} and \eqref{ibp4}, respectively.

  \begin{eqnarray}\label{ibp1}
    \int_{0}^{x}u(t) D_{q}v(t)\,d_{q}t &=&  \left[(uv)(t)\right]^{x}_{0}- \int_{0}^{x} D_{q} u(t) v(qt)\,d_{q}t,\\
\label{ibp2}
    \int_{0}^{x}u(x) D_{q^{-1}}v(x)\,d_{q} x&=&q \left[(uv)(t)\right]^{\frac{x}{q}}_{0} - \int_{0}^{x}D_{q^{-1}}u(t)v(\frac{t}{q})\,d_{q} t,
  \end{eqnarray}
  provided that   $(uv)(\cdot)$ is $q$-regular at zero. See also~\cite{annabyqq12}.

For $a>0,$ we use the notations $A_{q,a}$ and $A^{*}_{q,a}$ to donate

\begin{equation*}\label{aa1}
  A_{q,a}:=\{ a q^{n}: n\in \N_{0}\},\quad A^{*}_{q,a}:=A_{q,a}\cup \{0\}.
\end{equation*}

In the following, we consider the set $X$ of all functions defined on $[0,a]$ that satisfy 

\begin{equation*}
 \int_{0}^{a}|f(t)|^2 d_{q}t<\infty.
\end{equation*}

We define an equivalence relation on $X$ by  
\begin{equation*}
  f\sim g \ \text{if and only if }\  f(aq^k)=g(aq^k),\ k\in \N_{0}.
\end{equation*}
Let $ L^2_{q}[0,a]:=\left\{[f]: f\in X\right\},$ associated with the inner product
\begin{equation}\label{ip1}
  \langle [f], [g]\rangle=\int_{0}^{a} f(t) \overline{g(t)} d_{q} t, \quad f \in [f],g \in [g].
\end{equation}
One can verify that $L^2_{q}[0,a]$ associated with the well defined inner product \eqref{ip1} is a Hilbert space; see~\cite{annaby2005}. The space $C_{q}^{2}[0, a]$ is the subspace of $L_{q}^{2}[0, a]$ of all functions $y(\cdot)$ defined on $[0, a]$ such that $y(\cdot), D_{q} y(\cdot)$ are continuous at zero; see\cite{p2}. In a similar manner, we can define an equivalence relation on $X\times X$ defined by   
\begin{equation*}   f\sim g \ \text{if and only if }\  f(aq^k,aq^m)=g(aq^k,aq^m),\ k,m\in \N_{0}. 
\end{equation*}

Let \( L^2_{q}([0,a]\times[0,a]):=\left\{[f]: f\in X\times X\right\} \), associated with the inner product \begin{equation}\label{ip1}   \langle [f], [g]\rangle=\int_{0}^{a} \int_{0}^{a} f(x,t) \overline{g(x,t)} d_{q} x\, d_{q} t, \quad f \in [f], g \in [g]. \end{equation}

\subsection{The $q$-Sturm-Liouville problems}

 In \cite{annaby2005,annabyq12}, Annaby and Mansour introduced the following $q$-Sturm--Liouville problem
\begin{equation}\label{slp3}
 L_{q}(y)= -\dfrac{1}{q} D_{q^{-1}} D_{q} y(x)+v(x) y(x)=\lambda y(x), \quad 0 \leqslant x \leqslant a<\infty, \quad \lambda \in \mathbb{C},
\end{equation}
with the boundary conditions
\begin{equation}\label{slp4}
  U_{1}(y):=a_{11} y(0)+a_{12} D_{q^{-1}} y(0)=0,
\end{equation}
\begin{equation}\label{slp5}
  U_{2}(y):=a_{21} y(a)+a_{22} D_{q^{-1}} y(a)=0,
\end{equation}
where $v(\cdot)$ is continuous at zero and $\left\{a_{i j}\right\}, i, j \in\{1,2\}$ are arbitrary real numbers.  Recall that two  functions $\left\{y_{1}, y_{2}\right\}$ form a fundamental set of solutions of \eqref{slp3} if and only if their $q$-Wronskian, which  is defined as
 \begin{equation}\label{wr}
   W_{q}\left(y_{1}, y_{2}\right)(x):=y_{1}(x) D_{q} y_{2}(x)-y_{2}(x) D_{q} y_{1}(x), \quad x \in[0, a],
 \end{equation}
does not vanish at any point of $[0, a]$. Annaby and Mansour \cite{annaby2005,annabyq12} demonstrated the existence and uniqueness of solutions through the following theorem:
      \begin{theorem}\label{th22}
  Let $c_1,c_2$ be complex constants. Equation $\eqref{slp3}$ has a unique solution in $C_{q}^{2}[0, a]$ which satisfies
  \begin{equation}\label{nt2}
    \varphi(0, \lambda)=c_1, \quad D_{q^{-1}} \varphi(0, \lambda)=c_2, \quad \lambda \in \mathbb{C}.
  \end{equation}
Moreover, $\varphi(x, \lambda)$ is entire in $\lambda$ for all $x \in[0, a]$ and has the following form:
\begin{equation*}\label{}
\varphi(x, \lambda)= c_1\cos(sx;q)+c_2 \frac{\sin(sx;q)}{s}+\frac{q}{s} \int_0^x S_{q}(x,qt) v(q t) \varphi(q t, \lambda)\  \mathrm{d}_q t,
\end{equation*}
where $S_{q}(x,t)$ is defined as
\begin{equation}\label{sl7}
  S_{q}(x, t)=\sin (s x ; q) \cos (st ; q)-\cos (s x ; q) \sin (s t ; q), \quad x\in [0,a].
\end{equation}
\end{theorem}
 Moreover, they proved that the eigenvalues $\{\lambda_n\}_{n=0}^\infty$ of the problem \eqref{slp3}-\eqref{slp5} are real and the eigenfunctions $\{\varphi(\cdot,\lambda_n)\}_{n=0}^\infty$ associated with distinct eigenvalues are orthogonal in $L_ q ^2 [0, a]$ and hence the parserval's identity 
    \begin{equation}\label{qpi}
      \parallel f\parallel^2=\sum_{n=1}^{\infty}\dfrac{|<f,\varphi_{n}>|^2}{\parallel \varphi_{n}\parallel^2},\ 
     \end{equation}
     holds for $f \in L_ q ^2 [0, a]$, 
     see \cite[Lemma~4.3]{annaby2005}.
     
 In \cite{annabyqq12}, Annaby \etal \ studied the existence and uniqueness of solutions to the problem \eqref{slp3}-\eqref{slp5} when the interval $[0, a]$ is replaced with $[0, \infty)$. They developed the $q$-Titchmarsh-Weyl theory for singular $q$-Sturm-Liouville problems and proved that for any two functions $y$ and $z$ in $C_q^2(0, a),$ the Green's identity has the following form for $x\in [0,a]$
\begin{equation}\label{lin}
  \int_{0}^{x}\left(y L_q z-z L_q y\right)d_{q}t = W_{q^{-1}}\left(y, z\right)(0)- W_{q^{-1}}\left(y, z\right)(x),
\end{equation}
where $L_q$ is the $q$-difference operator defined in \eqref{slp3}, see~\cite{annabyqq12}.
\begin{remark}\label{eu1}
 Equation $\eqref{slp3}$ can be written as
     \begin{equation}\label{fff3}
        y(q x)+q y(\frac{x}{q})-\left[(1+q)+(v(x)-\lambda) x^2 (1-q)^2 \right]y(x)=0.
 \end{equation}
From Theorem~\ref{th22}, if $\eqref{slp3}$ is solved with the initial conditions $y(0)= D_{{q}^{-1}}y(0) = 0$, then the solution is the trivial solution using \eqref{fff3}. Also, if we solve equation $\eqref{slp3}$ with the boundary conditions $y(a) = D_{q^{-1}} y(a) = 0$, we conclude that $y(aq^ n) = y(0) = 0$ from the functional equation \eqref{fff3} and the continuity of the solution at zero.
\end{remark}
From Theorem~\ref{th22}, we deduce that if $ \varphi_0(x,\lambda)$ is the solution of 
\begin{equation}\label{pr1}
  -\dfrac{1}{q} D_{q^{-1}} D_{q} y(x)=\lambda y(x), \quad 0 \leqslant x \leqslant a,
\end{equation}
with the initial conditions \eqref{nt2} and $\varphi(x,\lambda)$ is the solution $\eqref{slp3}$ with the initial conditions \eqref{nt2}, then
$$
\varphi_0=(I+T) \varphi,
$$
where $I$ is the identity operator and $T$ is the operator defined as
\begin{eqnarray}\label{thi1}
 \nonumber T &:& L_q^2[0, a] \rightarrow L_q^2[0, a] \\
  && \varphi \rightarrow \int_0^x \dfrac{-S_{q}(x, t)}{s} v(t) \varphi(t, \lambda) d_{q} t,
\end{eqnarray}
in which $S_{q}(x, t)$ is defined in \eqref{sl7}. \\
In the next theorem, we prove that the operator $(I+T)$ is invertible.
\begin{theorem}\label{nt11s}
If $T$ is the operator defined in $\eqref{thi1}$, then $(I+T)$ is invertible and
\begin{equation*}\label{}
  (I+T)^{-1}=(I+E), 
\end{equation*}
where $E$ is the operator defined as
\begin{eqnarray}\label{thi2}
 \nonumber E &:& L_q^2[0, a] \rightarrow L_q^2[0, a] \\
 & & f \longrightarrow \int_0^x f(t) W(x, t)\ d_{q} t.
\end{eqnarray}
The kernel function $W(x, t)$ is defined by
  \begin{equation}\label{ks}
    W(x,t)=\sum_{n=1}^{\infty}W_{n}(x,t),
  \end{equation}
   where
    \begin{equation}\label{ks1}
      W_{1}(x,t)=\dfrac{1}{s} S_{q}(x,t) v(t),\quad  W_{n+1}(x, t)=\dfrac{q}{s} \int_{t}^{x} S_{q}(x, q\eta) v(q\eta) W_{n}(q\eta, t)\, d_{q}\eta,\quad n\in \N.
  \end{equation}
  In addition,
 \begin{equation}\label{nt3888}
   \|W\|^2=\int_{0}^{a} \int_{0}^{a} |W(x,t)|^2\, d_{q}t\, d_{q}x < \infty.
 \end{equation}
\end{theorem}
\begin{proof}
The general solution of~\eqref{pr1} with the initial conditions \eqref{nt2} is
\begin{equation}\label{nt7}
  \varphi_{0}(x, \lambda)= c_1\cos (s x; q)+\dfrac{c_2}{s} \sin (s x ; q).
\end{equation}
From Theorem~\ref{th22}, the general solution of \eqref{slp3} with the initial conditions \eqref{nt2} is given by
 \begin{equation}\label{nt10}
  \varphi(x, \lambda)=c_1 \cos (s x ; q)+\dfrac{c_2}{s} \sin (s x ; q)+\frac{q}{s} \int_{0}^{x} S_{q}(x, q t) v(q t) \varphi(q t, \lambda)\, d_{q} t,
\end{equation}
where $S_{q}(x,t)$ is defined as in \eqref{sl7}. Substitution with $\xi=q t$ in the $q$-integral~\eqref{nt10} and \eqref{nt7} yields
\begin{equation}\label{nt11}
 \varphi(x, \lambda)= \varphi_{0}(x, \lambda)+\frac{1}{s} \int_{0}^{qx} S_{q}(x,\xi) v(\xi)\varphi(\xi, \lambda)\, d_{q}\xi.
\end{equation}
Since $S_{q}(x,x)=0,$ from \eqref{sl7}, then \eqref{nt11} takes the form
\begin{align}\label{nt12}
     \varphi(x, \lambda):= & \varphi_{0}(x, \lambda)+\frac{1}{s} \int_{0}^{x} S_{q}(x,\xi) v(\xi) \varphi(\xi, \lambda)\, d_{q}\xi.
\end{align}
Using the method of successive approximations in \eqref{nt12}. Then, we have
\begin{equation}\label{nt13}
 \varphi_{0}(x, \lambda)=\varphi_{0}(x, \lambda),\ \quad \varphi_{n+1}(x, \lambda)=\frac{1}{s} \int_{0}^{x} S_{q}(x,\xi) v(\xi)\varphi_{n}(\xi, \lambda) \,d_{q}\xi \ (n\in\N_{0}).
\end{equation}
We shall prove the  identity 
\begin{equation}\label{nt14}
  \varphi_{n}(x, \lambda)=\int_{0}^{x} W_{n}(x, t) \varphi_{0}(t, \lambda) \,d_{q} t,\quad  W_{n}(x,t)\in L^2_q([0, a] \times[0, a]),
\end{equation}
where
\begin{equation*}
  \nonumber W_{n}(x, t)=\dfrac{1}{s} \int_{qt}^{qx} S_{q}(x, \xi) v(\xi) W_{n-1}(\xi, t)\, d_{q} \xi,\quad n>1,
\end{equation*}
 by using the mathematical induction. First, the basis of induction. Substituting $n=0$ in \eqref{nt13}, we obtain
 \begin{equation*}\label{nt15}
   \varphi_{1}(x, \lambda)=\frac{1}{s} \int_{0}^{x} S_{q}(x,t) v(t) \varphi_{0}(t, \lambda)\, d_{q}t.
\end{equation*}
 Thus, \eqref{nt14} is valid for $n=1$, where
\begin{equation}\label{nt16}
   W_{1}(x,t)=\dfrac{1}{s} S_{q}(x,t) v(t).
\end{equation}
Since $v(.)$  is continuous at zero and $S_{q}(x,t)$ is a continuous function for $(x,t)$ in $[0, a] \times[0, a]$, then $ W_{1}(x,t)\in L^2_q([0, a] \times[0, a]).$

Second, we suppose that \eqref{nt14} is valid for $n\geq1$ and $ W_{n}(x,t)\in L^2_q ([0, a] \times[0, a]).$

Third, we substitute \eqref{nt14} into \eqref{nt13} to obtain
\begin{equation}\label{nt17}
  \varphi_{n+1}\left(x, \lambda\right)=\frac{1}{s} \int_{0}^{x} S_{q}(x, \xi) v(\xi) \int_{0}^{\xi} W_{n}(\xi, t) \varphi_{0}(t, \lambda)\, d_{q}t \, d_{q} \xi.
\end{equation}
Since  $v(x),\ \varphi(x, \lambda)$ are $L^2_{q}[0, a]$ functions and $S_{q}(x, t),\  W_{n}(x, t) \in L^2_q([0, a] \times[0, a]),$ then the double $q$-integral \eqref{nt17} is absolutely convergent for all $(x,t)\in A_{q,a}\times A_{q,a}.$ Therefore, we can interchange the order of $q$-integrations to obtain
\begin{equation}\label{nt19}
 \varphi_{n+1}\left(x, \lambda\right)=\frac{1}{s} \int_{0}^{x}  \varphi_{0}(t, \lambda) \int_{qt}^{x} S_{q}(x, \xi) v(\xi)  W_{n}(\xi, t)\, d_{q}\xi \,d_{q}t.
\end{equation}
Therefore, by comparing \eqref{nt19} with \eqref{nt14}, we obtain
\begin{align}\label{nt20}
\nonumber W_{n+1}(x, t)=&\frac{1}{s} \int_{qt}^{x} S_{q}(x, \xi) v(\xi) W_{n}(\xi, t)\, d_{q} \xi\\
=&\frac{1}{s} \int_{qt}^{qx} S_{q}(x, \xi) v(\xi) W_{n}(\xi, t)\, d_{q} \xi.
\end{align}
Making the substitution $\xi=q\eta$ into the $q$-integral \eqref{nt20} yields
\begin{equation}\label{nt21}
W_{n+1}(x, t)=\dfrac{q}{s} \int_{t}^{x} S_{q}(x, q\eta) v(q\eta) W_{n}(q\eta, t) \,d_{q}\eta.
\end{equation}
Consequently, from the mathematical induction, \eqref{nt14} is true for all $n \in \N$ and $$ W_{n+1}(x,t)\in~L^2_q ([0, a]\times[0, a]).$$
Set 
\begin{equation}\label{su1}
  \varphi(x, \lambda):=\sum_{n=0}^{\infty}\varphi_{n}(x, \lambda).
\end{equation}
Then, from \eqref{ks}, we conclude 

\begin{equation*}
  \varphi(x, \lambda):=\varphi_{0}(x, \lambda)+\int_{0}^{x} W(x, t) \varphi_{0}(t, \lambda)\, d_{q} t,
\end{equation*}
 where $ W(x, t)$ is defined in \eqref{ks}.

Now, using the mathematical induction, we prove that
 \begin{equation}\label{nt29}
 | W_{n}(x,t)|\leq \dfrac{M^{n} L^{n}(1-q)^{n-1}q^{\frac{n(n-1)}{2}}}{|s|^{n}(q;q)_{n-1}} x^{n-1},\quad (x ,t) \in A_{q,a}\times A_{q,a},\ n\in \N,
\end{equation}
where $M$ and $L$ are the constants defined by
\begin{equation*}
  M=\max_{(x,t)\in ([0, a] \times[0, a])} |S_{q}(x, t)|,\ L=\max_{n\in \N_0} | v(aq^n)|.
\end{equation*}
$M$ exists because $S_{q}(x, t)$ is continuous on $[0, a] \times[0, a]$ and $L$ exists because $v(.)$ is continuous at zero. Using \eqref{nt16} for $t\leq x$, we have
\begin{equation*}
  \left|W_{1}(x, t)\right|=\dfrac{1}{|s|}| S_{q}(x, t)| | v(t)| \leq \dfrac{M L}{|s|}.
\end{equation*}
Therefore, $\eqref{nt29}$ is true at $n=1.$ We suppose that \eqref{nt29} holds for $n\geq1.$\\
From \eqref{nt21} for all $x ,t$ in $A_{q,a}$, we have
\begin{align}\label{nt28}
\nonumber| W_{n+1}(x,t)|\leq &\dfrac{q}{s} \int_{t}^{x}| S_{q}(x, q\eta)|| v(q\eta)| |W_{n}(q\eta, t)| \,d_{q}\eta\\
\leq&\dfrac{q M L}{|s|} \int_{t}^{x} | W_{n}(q\eta, t)| \,d_{q}\eta.
\end{align}
Using \eqref{nt29} in \eqref{nt28} for $t\leq x$, we have
\begin{align*}
\nonumber | W_{n+1}(x,t)|\leq&\dfrac{M^{n+1} L^{n+1}(1-q)^{(n-1)} q^{\frac{n(n-1)}{2}}}{|s|^{n+1} (q;q)_{n-1}} \int_{t}^{x} (q \eta)^{n-1} \,d_{q}\eta\\
 \leq &\dfrac{M^{n+1} L^{n+1} (1-q)^{n}q^{\frac{n(n+1)}{2}}}{|s|^{n+1} (q;q)_{n}} x^{n}.
\end{align*}
Hence, from the mathematical induction, \eqref{nt29} is correct for all $n \in \N.$
Consequently,
\begin{equation}\label{nt388}
  | W(x,t)| \leq \sum_{n=1}^{\infty}|W_{n}(x,t)| \leq \sum_{n=1}^{\infty}\dfrac{(M L)^{n} (1-q)^{(n-1)}q^{\frac{n(n-1)}{2}}}{|s|^{n} (q;q)_{n-1}} x^{n-1}= \dfrac{M L}{|s|}E_{q}(\frac{MLxq}{s}).
\end{equation}
From~\eqref{nt388}, we can obtain~\eqref{nt3888} and hence the series of \eqref{ks} converges absolutely and uniformly on $A_{q,a}$. Using \eqref{nt29} in \eqref{nt14}, we can prove that $\varphi(x,\lambda)$ is well defined  and 
\begin{equation}\label{nt3388}
  | \varphi(x,t)| \leq \sum_{n=1}^{\infty}|\varphi_{n}(x,t)| \leq \sum_{n=1}^{\infty}\dfrac{c_3(M L)^{n}(1-q)^n q^{\frac{n(n-1)}{2}}}{|s|^{n} (q;q)_n} x^{n}=c_3 E_{q}(\frac{MLx}{s})<\infty,
\end{equation}
where $|\varphi_{0}(x, \lambda)|<c_3$ and $c_3$ is a positive constant. Consequently, the series of \eqref{su1} converges absolutely and uniformly on $A_{q,a}.$ 
 \end{proof}

\section{A $q$-analog of Levinson-Marchenko uniqueness theorem}
 In this section, we introduce a $q$-analog of the Levinson-Marchenko uniqueness theorem, namely Theorem B. To achieve our goal, we need the following prerequisites.  We define a function $f$ to be equal to zero almost everywhere, and we denote it by $\left( f=0 \ (q a.e.) \right)$, if $f(x)$ is equal to zero only on a finite subset of $A_{q,a}$. We shall also use the notation $L_q(v,h,H)$ to donate the operator defined in \eqref{slp3} with the boundary conditions
 
\begin{equation}\label{is4}
  D_{q^{-1}} y(0)-h y(0)=0, \quad  D_{q^{-1}} y(a)+H y(a)=0.
\end{equation}
Moreover, we define $ \alpha_{n}$ as
\begin{equation}\label{is333}
      \alpha_{n}:=\int_{0}^{a}|\phi(x,\lambda_{n})|^2\, d_{q}x,\ n\in \N_0.
 \end{equation}
 where $\{\phi(x,\lambda_{n})\}_{n=0}^\infty$ are the eigenfunctions associated with the eigenvalues $\{\lambda_{n}\}_{n=0}^\infty$. The operator $\tilde{L}_q=L_q(\tilde{v},\tilde{h},\tilde{H})$ is defined similarly as the operator $L_q(v,h,H)$ as the following 
 \begin{equation}\label{sllp3}
\tilde{L}_q(y)= -\dfrac{1}{q} D_{q^{-1}} D_{q} y(x)+\tilde{v}(x) y(x)=\lambda y(x), \quad 0 \leqslant x \leqslant a, \quad \lambda \in \mathbb{C},
\end{equation}
 with replacing $v,h,H$ by $\tilde{v},\tilde{h},\tilde{H},$ respectively, and the boundary conditions are
 \begin{equation}\label{iss4}
  D_{q^{-1}} y(0)-\tilde{h} y(0)=0, \quad  D_{q^{-1}} y(a)+ \tilde{H} y(a)=0.
\end{equation}
 
  In the following, we discuss some properties of $q$-integral operators that assist us in proving our theorems.
 
\begin{theorem}\label{in1}
 Let $\chi: L^2_{q}[0, a] \rightarrow L^2_{q}[0, a]$ be the bounded linear operator defined by
 \begin{equation*}\label{in11}
   (\chi f)(t)=\int_0^a k(t, x) f(x)\, d_{q} x,
 \end{equation*}
where $k(t,x)\in L^2_q([0, a] \times[0, a])$. Then the adjoint operator $\chi^{*} $ is defined by
\begin{equation*}\label{}
  \left(\chi^{*} g\right)(t)=\int_0^a \overline{k(x, t)} g(x)\, d_{q}x,\quad g \in L^2_q[0, a].
\end{equation*}
\end{theorem}

\begin{proof}
The proof is similar to the classical case for the integral operator and is omitted. See~\cite[Theorem 11.2]{gohberg13}.
\end{proof}

\begin{theorem}\label{il1}
Let $A$ be the  bounded linear operator defined on $L^2_{q}[0, a]$ by
 \begin{equation}\label{nt77}
  A f(x)=f(x)+\int_{qx}^{a} Q(t, x) f(t)\, d_{q} t,\quad 0<x<t,
\end{equation}
 where $Q(t,x)$ is a real-valued function defined on $L^2_q([0, a] \times [0, a])$. Then the adjoint of the operator $A$ is given by
\begin{equation*}\label{ip66}
  A^*f(x)=f(x)+\int_0^x \overline{Q(x, t)} f(t)\, d_{q} t.
\end{equation*}
\end{theorem}

\begin{proof}
 Since $0<x<t,$ we may assume that $ Q(x, x)=0$ for all $x\in A_{q,a}$ and \eqref{nt77} is equivalent to
 \begin{equation}\label{ntt77}
  A f(x)=f(x)+\int_{x}^{a} Q(t, x) f(t)\, d_{q} t,\quad 0<x<t.
\end{equation}
Therefore, the proof follows from Theorem~\ref{in1} by taking
\begin{equation*}\label{}
  k(t, x):= \begin{cases}0, & t\leq x, \\ Q(t, x), & t > x.\end{cases}
\end{equation*}
\end{proof}
\begin{prop}\label{inv7}
   Let $A$ be the bounded operator defined on $L^2_{q}[0, a]$ by the following equation
   \begin{equation}\label{iil2}
     (Af)(x)=f(x)+\int_{x}^{a} Q(t, x) f(t)\,  d_{q} t, \quad x\in A_{q,a},
   \end{equation}
  in which $Q(t,x) \in L^2_q([0, a] \times[0, a])$. Then $A$ is an invertible operator and its inverse operator $A^{-1}$ is given by
 \begin{equation}\label{nt777}
(A^{-1} g)(x)=g(x)+\int_{x}^a R(t, x) g(t)\, d_q t,\quad g\in L^2_{q}[0, a],
\end{equation}
where $R(t,x)$ is defined recursively for $t, x$ in $A_{q,a}$ by
\begin{eqnarray}\label{valr}
 \nonumber  R(t,a)&=& 0, \\
  R(t,aq^m) &=& - Q(t, a q^m)-\int_{aq^m}^{t} Q(s, x) R(t, s)\,  d_{q} s,
\end{eqnarray}
and the $q$-integral in \eqref{valr} is considered zero if $t\in\{aq^k, k\leq m\}.$ 
 \end{prop}
\begin{proof}
Set $Af=g.$ Then from \eqref{iil2}, we have
  \begin{equation}\label{nt6633}
        f(x)=g(x)-\int_{x}^{a} Q(t, x) f(t)\,  d_{q} t.
  \end{equation}
  Substituting $x=a$ into \eqref{nt6633} yields
  \begin{equation}\label{nt45}
    f(a)=g(a).
  \end{equation}
  We shall prove using the strong mathematical induction that 
  \begin{equation}\label{il33}
       f(aq^j)=g(aq^j)+\int_{aq^j}^{a} R(t, aq^j) g(t)\,  d_{q} t,\ j\in \mathbb{N}_0.
  \end{equation}
  First, the basis of induction, from equation~\eqref{nt45}, \eqref{il33} holds at $j=0$. Second, we assume that~\eqref{il33} holds for $0\leq j\leq m.$ Hence, from \eqref{nt6633}, we have
  \small
  \begin{align}\label{il433}
     \nonumber  f(aq^{m+1})=& g(aq^{m+1})-\int_{aq^{m+1}}^{a} Q(s, a q^{m+1}) f(s)\,  d_{q} s\\
      \nonumber =& g(aq^{m+1})-\int_{aq^{m+1}}^{a} Q(s, a q^{m+1})\left(g(s)+\int_{s}^{a} R(t, s) g(t)\,  d_{q} t\right)\,  d_{q} s\\
       =& g(aq^{m+1})-\int_{aq^{m+1}}^{a} Q(t, a q^{m+1}) g(t)\,  d_{q} s-\int_{aq^{m+1}}^{a} Q(s, a q^{m+1})\int_{s}^{a} R(t, s) g(t)\,  d_{q} t\,  d_{q}s.
       \end{align}
  Since  $g(t)$ is a $L^2_{q}[0, a]$ function and $Q(x,t)\in L^2_q([0, a] \times[0, a]),$ then the double $q$-integral on the right hand side of \eqref{il433} is absolutely convergent. Therefore, we can interchange the order of $q$-integration to obtain
  \begin{eqnarray*}
     f(aq^{m+1}) &=& g(aq^{m+1})+\int_{aq^{m+1}}^{a}g(t) \left(-Q(t, a q^{m+1})- \int_{aq^{m+1}}^{t} Q(s, a q^{m+1})R(t, s)\,  d_{q} s\right)\,  d_{q}t\\
     &=& g(aq^{m+1})+\int_{aq^{m+1}}^{a} g(t) R(t, aq^{m+1})\,  d_{q}t.
  \end{eqnarray*}

Hence, from the mathematical induction, the identity~\eqref{il33} holds for all $j\in\mathbb{N}_0$ and \eqref{valr} is obtained.  
\end{proof}
\begin{prop}\label{thsao1}
 Let $T$ be a self-adjoint operator on a complex vector inner product space $V$. If $\langle T v, v\rangle=0$ for all $v \in V$, then $T=0$.
\end{prop}
\begin{proof}
  The proof follows from the fact that if $T$ is a self adjoint operator on a complex vector inner product space $V$, then 
  \begin{equation*}
    \|T\|=\sup_{v\in V}|\langle T v, v\rangle|,
  \end{equation*}
 see~\cite[Corollary 4.2]{gohberg13}.
\end{proof}

 Now we introduce a $q$-analog of Theorem B, the Levinson-Marchenko theorem.
  
 \begin{theorem}\label{nt60} 
 Let $\{\lambda_{n}, \alpha_{n} \}_{n=0}^\infty$ and $\{\tilde{\lambda}_{n}, \tilde{\alpha}_{n}\}_{n=0}^\infty$ be the spectral data of the problems
 $L_q(v,h,H)$ and $L_q(\tilde{v},\tilde{h},\tilde{H})$, respectively. If $\lambda_{n}=\tilde{\lambda}_{n}, \alpha_{n}=\tilde{\alpha}_{n}, n \in \N_{0}$, and $\tilde{h}=h,$ then $\tilde{L}_q=L_q$, i.e., $\tilde{v}(x)=v(x)$ ($q$a.e.) and $\tilde{H}=H$.
\end{theorem}

\begin{proof}
If $\varphi$ is a solution of $L_q(v,h,H),$ then from Theorem~\ref{th22} and Theorem~\ref{nt11s}, we have
\begin{equation}\label{thh21}
  \varphi(x, \lambda)=(I+E) \psi_0(x, \lambda), \quad \psi_0(x, \lambda)=(I+T) \varphi(x, \lambda),
\end{equation}
where $T$ and $E$ are the operators defined in \eqref{thi1} and \eqref{thi2}, respectively. Here, the function $\psi_{0}$ is the solution of \eqref{pr1} with the initial conditions $\psi_{0}(0)=1, D_{q^{-1}}\psi_{0}(0)=h$.\\
Similarly if $\tilde{\varphi}$ is a solution of $L_q(\tilde{v},\tilde{h},\tilde{H})$ with $\tilde{h}=h,$ then 
\begin{equation}\label{thh22}
 \tilde{\varphi}(x, \lambda)=(I+\tilde{E}) \psi_0(x, \lambda), \quad \psi_0(x, \lambda)=(I+\tilde{T}) \varphi(x, \lambda),
\end{equation}
where $\tilde{E}$ and $\tilde{T}$ are the operators defined on $L_q^2[0,a]$ as
\begin{eqnarray}
  \nonumber  (\tilde{E}{f})(x) &:=& \int_{0}^{x} \tilde{W}(x, t){f}(t)\, d_{q} t, \\
 \nonumber (\tilde{T}{f})(x)&:=&\int_0^x \dfrac{-S_{q}(x, t)}{s} \tilde{v}(t) {f}(t)\, d_{q} t.
\end{eqnarray}

From \eqref{thh21} and \eqref{thh22}, we have 
\begin{align}\label{nt666}
  \nonumber  \varphi(x, \lambda)=&\left(I+E \right) (I+\tilde{T}) \tilde{\varphi}(x, \lambda)\\
    =&(I+E+\tilde{T}+E\tilde{T})\tilde{\varphi}(x, \lambda).
\end{align}
Equation \eqref{nt666} becomes
\begin{equation}\label{nt68}
  \varphi(x, \lambda)=\tilde{\varphi}(x, \lambda)+ \int_{0}^{x} Q(x,t) \tilde{\varphi}(t, \lambda)\, d_{q} t,
\end{equation}
where $Q(x,t)$ is the real-valued function defined for $(x,t)\in A_{q,a}\times A_{q,a}$ by
\begin{equation*}
   Q(x,t)= W(x, t)-\dfrac{S_{q}(x, t)}{s} \tilde{v}(t)-\frac{1}{s}\int_{qt}^{x} W(x,\xi)S_{q}(\xi, t) \tilde{v}(t)\, d_{q} \xi.
\end{equation*}
From~\eqref{sl7}, \eqref{ks} and \eqref{ks1}, we conclude that  $$W(x,x)=S_{q}(x, x)=0, \quad ( x \in A_{q,a}).$$ Therefore,
\begin{equation*}\label{QQ}
  Q(x,x)=0, \quad (x \in A_{q,a}).
\end{equation*}
Let $f \in L^2_{q}[0, a]$. From \eqref{nt68}, we have
\begin{equation}\label{nt70}
  \int_{0}^{a} f(x) \varphi(x, \lambda) d_{q} x= \int_{0}^{a} f(x) \tilde{\varphi}(x,\lambda)\, d_{q} x+\int_{0}^{a} f(x) \int_{0}^{x} Q(x,t) \tilde{\varphi}(t, \lambda)\, d_{q} t\, d_{q} x.
\end{equation}
Since  $f(x),\ \tilde{\varphi}(x, \lambda)$ are $L^2_{q}[0, a]$ functions and $Q(x,t)\in L^2_q([0, a] \times[0, a]),$ then the double $q$-integral on the right hand side of \eqref{nt70} is absolutely convergent. Therefore, we can interchange the order of $q$-integration to obtain
 \begin{equation}\label{nt71}
  \int_{0}^{a} f(x) \varphi(x, \lambda) d_{q} x= \int_{0}^{a} f(x) \tilde{\varphi}(x,\lambda)\, d_{q} x+\int_{0}^{a}  \int_{qx}^{a}f(t) Q(t,x) \tilde{\varphi}(x, \lambda)\, d_{q} t\, d_{q} x.
\end{equation}
Equation~\eqref{nt71} can be written as
\begin{equation}\label{nt72}
  \int_{0}^{a} f(x) \varphi(x, \lambda)\, d_{q} x=\int_{0}^{a} g(x) \tilde{\varphi}(x, \lambda)\, d_{q} x,
\end{equation}
where $g(x)=f(x)+\int_{qx}^{a} Q(t, x) f(t)\, d_{q} t.$ Set
\begin{gather}\label{nt74}
b_{n}:=\int_{0}^{a} f(x) \varphi\left(x, \lambda_{n}\right)\, d_{q} x, \quad d_{n}:=\int_{0}^{a} g(x) \tilde{\varphi}\left(x, \lambda_{n}\right)\, d_{q} x,\quad n\in \N_{0}.
\end{gather}
From \eqref{nt72} and \eqref{nt74}, we obtain $b_{n}=d_{n}$ for all $n\in \N_{0}.$ Using Parseval's equality \eqref{qpi} and \eqref{nt74}, we conclude that
\begin{equation}\label{nt75}
  \int_{0}^{a}|f(x)|^{2}\, d_q x=\sum_{n=0}^{\infty} \dfrac{\left|\int_{0}^{a} f(x) \varphi\left(x, \lambda_{n}\right)\, d_{q} x\right|^{2}}{\int_{0}^{a} \varphi^{2}\left(x, \lambda_{n}\right)\, d_{q} x}=\sum_{n=0}^{\infty} \dfrac{\left|b_{n}\right|^{2}}{\alpha_{n}}=\sum_{n=0}^{\infty} \dfrac{\left|{d}_{n}\right|^{2}}{\tilde{\alpha}_{n}}=\int_{0}^{a}|g(x)|^{2}\, d_{q} x.
\end{equation}
Then, from \eqref{nt75}, we have
\begin{equation}\label{nt76}
\|f\|^2=\|g\|^2.
\end{equation}
Since $A f(x)=g(x),$ then $\|A f\|^2=\|f\|^2$ for any $f\in L^2_{q}[0, a].$ Consequently
   \begin{equation}\label{io1}
       \langle A f, Af\rangle=\left\langle f, f\right\rangle.
   \end{equation}
   Equation~\eqref{io1} becomes
\begin{equation*}\label{io2}
       \langle (A^* A-I) f,f\rangle=0, \text{for all}\  f\in L_{q}^{2}[0, a].
   \end{equation*}
    Since the operator $T:=(A^* A-I)$ is a self-adjoint operator on the complex inner product space $L_{q}^{2}[0, a]$, then from Proposition~\ref{thsao1}, we have $A^* A=I$ and hence $A^{-1}=A^*.$ Therefore, using Theorem~\ref{il1} and Proposition~\ref{inv7}, we conclude that
\begin{equation}\label{pp4}
  \int_0^x \overline{Q(x, t)} g(t)\, d_{q} t=\int_{ x}^a R(t, x) g(t)\, d _{q} t \quad \text { for all }\ g \in L_{q}^2[0,a],\ x\in A_{q,a},
\end{equation}
where $ R(t, x)$ is defined in \eqref{valr}. For $k\in \N,$ define the function $g_{k}$ by
\begin{equation*}\label{pp5}
  g_{k}(t)=\left\{\begin{array}{c}0,\quad  t \neq a q^k, \\ 1,\quad t=a q^k.\end{array}\right.
\end{equation*}
In \eqref{pp4} set $x=aq^m,\ m\in \N$ and if $k \geq m$, substitute with $g=g_{k}$ on \eqref{pp4}, then  
\begin{equation*}
  a q^{k+m}(1-q) \overline{Q\left(a q^m, aq^{k}\right)}=0.
\end{equation*}
Hence,
\begin{equation*}
   Q\left(a q^m, aq^{k}\right)=0 \ \text{for all}\  k \geq m.
\end{equation*}
Therefore, from \eqref{nt68}, we conclude that
\begin{equation}\label{dd1}
  \varphi(x, \lambda) = \tilde{\varphi}(x, \lambda),\quad x\in A_{q,a}.
\end{equation}
Consequently, using \eqref{slp3} and \eqref{sllp3}, we can obtain 
\begin{equation*}\label{}
  (v(x)-\tilde{v}(x))\varphi(x,\lambda)=0,\ \text{for all} \ x\in A_{q,a}.
\end{equation*}
Since $\varphi(x,\lambda)$ is continuous at zero and $\varphi(0,\lambda)\neq 0,$ then $\varphi(x,\lambda)=0$ possibly at a finite subset of $A_{q,a}.$ Therefore, $v(x)=\tilde{v}(x)$ (qa.e.). In addition, using \eqref{is4}, \eqref{iss4} and \eqref{dd1}, we have $\left(H-\tilde{H}\right)\varphi(a,\lambda)=0.$ From Remark~\ref{eu1}, we have $\varphi(a)\neq0.$ Hence $H=\tilde{H}.$
\end{proof}
\section{A $q$-analog of the Gelfand-Levitan equation with an application}
 
In this section, we introduce a q-analog of the Gelfand-Levitan equation and use it to transform the $q$-difference operators into $q$-integral operators; see~\cite{strum,gel51,freiling01}. This helps us to prove a $q$-analog of the Ashrafyan uniqueness theorem, namely Theorem C.

The Gelfand-Levitan equation is defined by
\begin{equation}\label{gle}
  G(x, t)+F(x, t)+\int_0^x G(x, s) F(s, t)\ d s=0, \quad 0<t<x,
\end{equation}
where $F(x, t)$ is called the kernel function. This equation was introduced by Gelfand and Levitan in \cite{gel51}. The Gelfand-Levitan method is a powerful tool for solving the inverse Sturm-Liouville problems. This method involves reconstructing the differential equation from its spectral data.
In \cite{ashrafyan17, jodeit97,freiling01}, the authors used the Gelfand-Levitan integral equation~\eqref{gle} to prove uniqueness theorems for inverse Sturm-Liouville problems.
 \vskip .5 cm
\par Let $v_{0}(x)$ be a continuous function  at zero defined on $[0,a]$, and $h_{0}$ and $H_{0}$ be fixed real numbers. Consider the $q$-Sturm-Liouville problem
\begin{equation}\label{p1}
   L_{q,0}:=L_q(v_{0},h_{0},H_{0})=-\dfrac{1}{q} D_{q^{-1}} D_{q} y(x)+v_{0}(x) y(x)=\lambda y,\quad  x \in [0,a],
\end{equation}
with the boundary conditions
  \begin{equation}\label{pp1}
  D_{q^{-1}} y(0)-h_{0} y(0)=0, \quad  D_{q^{-1}} y(a)+H_{0} y(a)=0.
\end{equation}
 Let $\phi_{0}(x, \lambda)$ be the solution of \eqref{p1} with the initial conditions
\begin{equation}\label{ppp1}
  \phi_{0}(0, \lambda)=1,\quad D_{q^{-1}}\phi_{0}(0, \lambda)=h_{0}.
\end{equation}
Then $\phi_{0}(x, \lambda)$ satisfies the first boundary condition in \eqref{pp1} for every $\lambda\in \C$.

The eigenvalues, $\{\lambda_{n,0}\}_{n=0}^\infty,$ of the problem $\eqref{p1}-\eqref{ppp1}$ are the roots of the equation
\begin{equation}\label{pppp1}
   D_{q^{-1}}\phi_{0}(a, \lambda)+H_{0} \phi_{0}(a, \lambda)=0.
\end{equation}
The corresponding eigenfunctions are $\{\phi_{0}\left(\cdot, \lambda_{n,0}\right)\}_{n=0}^\infty$ and the normalization constants $\left\{\alpha_{n,0}\right\}$ are defined by
\begin{equation}\label{nc1}
  \alpha_{n, 0}:=\int_{0}^{a} |\phi_{0}\left(x, \lambda_{n,0}\right)|^2\, d_{q} x,\quad n\in \N_{0}.
\end{equation}
The $q$-Gelfand-Levitan equation has the following form
\begin{equation}\label{i2}
  K(x, y)+F(x, y)+\int_{0}^{x} K(x, t) F(t, y)\, d_{q} t=0, \quad (x,y)\in A_{q,a}\times A_{q,a},\  0<y<x.
\end{equation}
 Let $\{c_{n}\}_{n=0}^\infty$ be a sequence of real numbers that is not identically zero and for which the series $\sum_{n=0}^{\infty} c_{n}$ converges. Define the kernel function $F(x, y)$  by
\begin{equation}\label{i1}
F(x, y):=\sum_{n=0}^{\infty} c_{n} \phi_{0}\left(x, \lambda_{n,0}\right) \phi_{0}\left(y, \lambda_{n,0}\right),\ (x,y)\in A_{q,q}\times A_{q,a}.
\end{equation}
\begin{theorem}\label{n1}
Let $\left\{\lambda_{n,0}, \alpha_{n,0}\right\}_{n=0}^\infty$ be the spectral date of the operator $L_{q,0}$ defined by $\eqref{p1}$. Assume that the coefficients $\{c_{n}\}_{n=0}^\infty$ of the function $F(x, y)$ defined by $\eqref{i1}$ satisfy
\begin{equation}\label{i3}
  \lim _{n \rightarrow \infty} c_n=0,\quad 1+c_{n} \alpha_{n, 0}>0,\ \text{where} \  \alpha_{n, 0} \ \text{is defined by \eqref{nc1}}.
 \end{equation}
 Then the $q$-Gelfand-Levitan equation $\eqref{i2}$ has a unique solution on $L_q^2([0,a]\times [0,a])$.
\end{theorem}
\begin{proof}
If $\eqref{i2}$ has two distinct solutions, then the homogenous equation
\begin{equation}\label{ii2}
  K(x, y)+\int_{0}^{x} K(x, t) F(t, y)\, d_{q} t=0,\quad 0<y<x,\ (x,y)\in A_{q,a}\times A_{q,a},
\end{equation}
has a non-zero solution. Therefore, to prove the uniqueness of the solution of $\eqref{i2}$, we prove that equation~\eqref{ii2} has only the trivial solution. Fix $x \in A_{q, a}$. Suppose that $h_{x} (y) \left( y \in A_{q, a}\right)$  is a nontrivial solution of \eqref{ii2}. That is, 
\begin{equation}\label{th1}
  h_{x}(y)+\int_{0}^{x} F(t, y) h_{x}(t)\, d_{q}t=0.
\end{equation}
 By multiplying equation~\eqref{th1} with $h_{x}(y)$ and $q$-integrating in $y$ from 0 to $x$, we obtain
\begin{equation}\label{th2}
  \int_{0}^{x} h_{x}^{2}(y)\, d_{q} y+\int_{0}^{x} \int_{0}^{x} F(t, y) h_{x}(t) h_{x}(y)\, d_{q}t \,d_{q} y=0.
\end{equation}
From \eqref{ii2}, $y<x$. Therefore, we may assume without any loss of generality that $h_{x} (y)=0$ for $y\geq x$. Hence
\begin{equation}\label{th222}
   \int_{0}^{x} h_{x}^{2}(y)\, d_{q} y= \int_{0}^{a} h_{x}^{2}(y)\, d_{q} y.
\end{equation}
Hence, from Parseval equality \eqref{qpi}, we obtain
\begin{equation}\label{th11}
   \int_{0}^{a} h_{x}^{2}(y)\, d_{q} y=\sum_{n=0}^{\infty} \frac{1}{\alpha_{n, 0}}\left(\int_{0}^{a} h_{x}(y) \phi_{0}\left(y, \lambda_{n,0}\right)\, d_{q} y\right)^{2}.
\end{equation}
Substituting \eqref{i1}, \eqref{th222}, \eqref{th11} into \eqref{th2} and using $h_{x} (y)=0$ for $y>x$, we obtain
\begin{align}\label{f5}
\nonumber&\int_{0}^{a} h_{x}^{2}(y) d_{q}y+\int_{0}^{a} \int_{0}^{a} F(t, y) h_{x}(t) h_{x}(y) d_{q}t d_{q}y \\
&\quad=\sum_{n=0}^{\infty} \frac{1}{\alpha_{n, 0}}\left[1+c_{n} \alpha_{n, 0}\right]\left[\int_{0}^{a} h_{x}(y) \phi_{0}\left(y, \lambda_{n,0}\right) d_{q}y\right]^{2}=0.
\end{align}
Using~\eqref{f5} and condition~\eqref{i3}, it follows that
\begin{equation}\label{ff5}
  \int_{0}^{a} h_{x}(y) \varphi_{0}\left(y, \lambda_{n,0}\right) d_{q}y=0, \quad n\in \N_{0}.
\end{equation}
From \eqref{ff5}, the function $h_{x}(y)$ is orthogonal to the complete orthogonal basis $\{\phi_{0}\left(y, \lambda_{n,0}\right)\}$. Hence $h_{x}(y)=0$ for all $y\in A_{q,a}$. Since $x \in A_{q,a}$ is arbitrary, then $h_{x}(y)=h(x,y)=0$ for all $(x,\ y) \in A_{q,a}\times A_{q,a}.$
\end{proof}
\begin{theorem}\label{T2}
 The solution of the $q$-Gelfand-Levitan equation~\eqref{i2} satisfies the following $q$-difference equation
\begin{equation}\label{l6}
   D_{q, x} D_{q^{-1}, x} K(x, y)- v(x) K(x, y)=  D_{q, y} D_{q^{-1}, y} K(x, y)- v_{0}(y) K(x, y),\quad (x,y)\in A_{q,a}\times A_{q,a},
\end{equation}
where
\begin{equation}\label{l7}
 v(x):= v_{0}(x)+D_{q, x}K(x, t)\mid_{t=x} + D_{q^{-1}, t} K(x, t)\mid_{t=x},\quad x\in A_{q,a},
\end{equation}
\begin{equation}\label{l8}
   h_{0} K(x, 0)- D_{q^{-1}, t} K(x, t)|_{t=0} =0,
\end{equation}
and
\begin{equation}\label{l99}
   K(x, x)=0,\quad \text{for all} \quad x \in A_{q,a}.
\end{equation}
\end{theorem}
\begin{proof}
Set
\begin{equation}\label{l1}
 J:=K(x, y)+F(x, y)+\int_{0}^{x} K(x,t) F(t, y) d_{q}{t},\quad (x,y)\in  A_{q,a}\times A_{q,a}.
\end{equation}
Applying $D_{q^{-1}, x}$ on~\eqref{l1} and using \eqref{ff1}, we obtain
\begin{align}\label{ll22}
  \nonumber D_{q^{-1}, x} J &=  D_{q^{-1}, x} K(x, y)+D_{{q}^{-1}, x} F(x, y) + \int_{0}^{x} D_{q^{-1}, x} K(x, t) {F(t, y)}\, d_{q} t  \\
     +& K\left(\frac{x}{q}, \frac{x}{q}\right) F\left(\frac{x}{q}, y\right).
\end{align}
Operating with $D_{q, x}$ on \eqref{ll22} and using \eqref{ff2} yields
\begin{eqnarray}\label{l3}
 \nonumber D_{q, x} D_{q^{-1}, x} J &=& D_{q, x} D_{q^{-1}, x} K(x, y)+ D_{q, x} D_{{q}^{-1}, x} F(x, y)+D_{q, x}  K( x, t)\mid_{t=x} {F(x, y)}\\
   &+& D_{q, x}\left[K\left(\frac{x}{q}, \frac{x}{q}\right) F\left(\frac{x}{q}, y\right) \right]+
\int_{0}^{x} D_{q, x} D_{q^{-1}, x} K(x, t) {F(t, y)}\, d_{q} t.
\end{eqnarray}
Similarly, applying $D_{q, y} D_{q^{-1}, y}$ on~\eqref{l1} yields
\begin{equation}\label{l4}
D_{q, y} D_{q^{-1}, y} J= D_{q, y}D_{q^{-1}, y} K(x, y)+ D_{q, y}D_{{q}^{-1}, y} F(x, y)+
\int_{0}^{x} K(x, t) D_{q, y} D_{q^{-1}, y} {F(t, y)}\, d_{q} t.
\end{equation}
From \eqref{i1} and \eqref{p1}, we obtain
\begin{align}\label{l5}
  D_{q, y} D_{q^{-1}, y} {F(t, y)}-D_{q, t} D_{q^{-1}, t} {F(t, y)}
  = & (v_{0}(y)-v_{0}(t)){F(t, y)}.
\end{align}
Substituting~\eqref{l5} into \eqref{l4}, we get
\begin{align}\label{ll6}
    D_{q, y} D_{q^{-1}, y} J=&  D_{q, y}D_{q^{-1}, y} K(x, y)+ D_{q, y} D_{{q}^{-1}, y} F(x, y) +I_{1}+I_{2},
\end{align}
where
\begin{align*}\label{}
 I_{1}:&=\int_{0}^{x} K(x, t) D_{q, t} D_{q^{-1}, t} {F(t, y)}\, d_{q} t,\\
 I_{2}:&=\int_{0}^{x}  (v_{0}(y)-v_{0}(t)) K(x, t) {F(t, y)}\, d_{q} t.
\end{align*}
From~\eqref{ibp1}, applying twice the $q$-integration by parts rule to $I_{1}$ yields
\begin{equation}\label{ll66}
I_{1}= [K(x, t) D_{q^{-1}, t} F(t, y)]^{x}_{t=0}-[ D_{q^{-1}, t} K(x, t) F(t, y)]^{x}_{t=0}+\int_{0}^{x}  D_{q, t} D_{q^{-1}, t} K(x, t) F(t, y)\, d_{q} t.
 \end{equation}
Substituting  \eqref{ll66} into \eqref{ll6}, we get
\begin{align}\label{ll8}
   \nonumber D_{q, y} D_{q^{-1}, y} J=&  D_{q, y}D_{q^{-1}, y} K(x, y)+ D_{q, y}D_{{q}^{-1}, y} F(x, y)\\
  \nonumber +&\int_{0}^{x}  \left[v_{0}(y)K(x, t)-v_{0}(t)K(x, t)+ D_{q, t} D_{q^{-1}, t}K(x, t)\right] {F(t, y)}\, d_{q} t \\
    +& [K(x, t) D_{q^{-1}, t} F(t, y)]^{x}_{t=0}-[ D_{q^{-1}, t} K(x, t) F(t, y)]^{x}_{t=0}.
\end{align}
 We obtain from equations \eqref{l1},\ \eqref{l3},\ \eqref{l5}, and \eqref{ll8} the following formula:
\begin{align*}
\nonumber & D_{q, x} D_{q^{-1}, x} J - D_{q, y} D_{q^{-1}, y} J+ v_{0}(y) J-v(x) J=\\
\nonumber&D_{q, x}D_{q^{-1}, x} K(x, y)- D_{q, y}D_{q^{-1}, y} K(x, y)+(v_{0}(y)-v(x))K(x, y)\\
\nonumber&+\left(v_{0}(x)+D_{q, x} K(x, t)\mid_{t=x}+D_{q^{-1}, t} K(x, t)\mid_{t=x}-v(x) \right) F(x, y)- K(x, x) D_{q^{-1}, t} F(t,y)|_{t=x}\\
\nonumber& +K(x, 0) D_{q^{-1}, t} F(t, y)|_{t=0}-D_{q^{-1}, t} K(x, t)|_{t=0} F(0, y)+ D_{q, x}\left[K\left(\frac{x}{q}, \frac{x}{q}\right) F\left(\frac{x}{q}, y\right) \right]\\
&+ \int_{0}^{x} \left [ D_{q, x} D_{q^{-1}, x} K(x, t)- D_{q, t} D_{q^{-1}, t} K(x, t)+ v_{0}(t) K(x, t)-v(x) K(x, t)\right] F(t, y)\, d_{q} t=0.
\end{align*}
By assuming \eqref{l7} and \eqref{l99}, we obtain
\begin{align}\label{ll999}
\nonumber&D_{q, x}D_{q^{-1}, x} K(x, y)- D_{q, y}D_{q^{-1}, y} K(x, y)+(v_{0}(y)-v(x))K(x, y)\\
\nonumber& +K(x, 0) D_{q^{-1}, t} F(t, y)|_{t=0}-D_{q^{-1}, t} K(x, t)|_{t=0} F(0, y)\\
&+ \int_{0}^{x} \left [ D_{q, x} D_{q^{-1}, x} K(x, t)- D_{q, t} D_{q^{-1}, t} K(x, t)+ v_{0}(t) K(x, t)-v(x)K(x, t)\right] F(t, y)\, d_{q} t=0.
\end{align}
 From \eqref{ppp1} and \eqref{i1}, we conclude that
\begin{equation}\label{lll3}
  D_{q^{-1}, t} F(t, y)|_{t=0}=h_{0}F(0, y),
\end{equation}
where $F(0, y)\neq0.$ Hence, from \eqref{lll3} and \eqref{l8}, we obtain 
\begin{equation}\label{pp11}
  K(x, 0) D_{q^{-1}, t} F(t, y)|_{t=0}-D_{q^{-1}, t} K(x, t)|_{t=0} F(0, y)=0.
\end{equation}
Substituting \eqref{pp11} into \eqref{ll999} yields
\begin{align*}
&D_{q, x}D_{q^{-1}, x} K(x, y)- D_{q, y}D_{q^{-1}, y} K(x, y)+(v_{0}(y)-v(x))K(x, y)\\
&+ \int_{0}^{x} \left [ D_{q, x} D_{q^{-1}, x} K(x, t)- D_{q, t} D_{q^{-1}, t} K(x, t)+ v_{0}(t) K(x, t)-v(x)K(x, t)\right] F(t, y)\, d_{q} t=0.
\end{align*}
Hence, the function
\begin{equation*}
  h_{x}(y)=D_{q, x}D_{q^{-1}, x} K(x, y)- D_{q, y}D_{q^{-1}, y} K(x, y)+(v_{0}(y)-v(x))K(x, y),
\end{equation*}
 is a solution of \eqref{th1}. Using Theorem~\ref{n1}, we conclude that $ h_{x}(y)=0$ for all $y \in A_{q,a}.$
This results in \eqref{l6} and concludes the proof.
\end{proof}

\begin{theorem}\label{nt}
  If $K(x, y)$ is the solution of the $q$-Gelfand-Levitan equation~\eqref{i2} and $\phi_{0}(x, \lambda)$ is a solution of~\eqref{p1}-\eqref{ppp1}, then for every complex number $\lambda,$ the function
  \begin{equation}\label{t7}
    \phi(x, \lambda):=\phi_{0}(x, \lambda)+\int_{0}^{x} K(x, t) \phi_{0}(t, \lambda)\, d_{q} t,\quad x\in A^{*}_{q,a},
  \end{equation}
is a solution of the $q$-difference equation~\eqref{slp3} subject to the initial conditions
\begin{equation}\label{t78}
  \phi(0, \lambda)=1, \quad D_{q^{-1}}\phi(0, \lambda)=h_{0}- q\sum_{k=0}^{\infty} c_{k}=h,
\end{equation}
with \begin{equation*}
 v(x):= v_{0}(x)+D_{q, x}K(x, t)\mid_{t=x} + D_{q^{-1}, t} K(x, t)\mid_{t=x},\quad x\in A_{q,a},
\end{equation*}
 \end{theorem}
 
 \begin{proof}
 Since $\phi_{0}(x, \lambda)$ is a solution of~\eqref{p1}-\eqref{ppp1}, then
\begin{align}\label{t3}
  \nonumber \lambda \phi(x, \lambda)=& \lambda \phi_{0}(x, \lambda)+\int_{0}^{x} \lambda \phi_{0}(t, \lambda) K(x, t)\, d_{q} t\\
  =&\lambda \phi_{0}(x, \lambda)-\int_{0}^{x}D_{q, t} D_{q^{-1}, t}\phi_{0}(t, \lambda)K(x, t)\, d_{q}t +\int_{0}^{x}v_{0}(t)\phi_{0}(t, \lambda) K(x, t)\, d_{q}t.
\end{align}
Using the $q$-integration by parts rule~\eqref{ibp1} yields
\begin{align*}\label{tt44}
\int_{0}^{x} K(x, t) D_{q, t} D_{q^{-1}, t} { \phi_{0}(t, \lambda)}\,  d_{q} t=& [K(x, t) D_{q^{-1}, t}  \phi_{0}(t, \lambda)]_{t=0}^x-[ D_{q^{-1}, t} K(x, t)  \phi_{0}(t, \lambda)]^{x}_{t=0}\\
&+\int_{0}^{x}  D_{q, t} D_{q^{-1}, t} K(x, t)  \phi_{0}(t, \lambda)\,  d_{q} t.
 \end{align*}
Hence, \eqref{t3} takes the form
\begin{eqnarray}\label{t5}
\nonumber  \lambda \phi(x, \lambda) &=& \lambda \phi_{0}(x, \lambda)+ \int_{0}^{x}\left[v_{0}(t)K(x, t)-D_{q, t} D_{q^{-1}, t} K(x, t)\right]\phi_{0}(t, \lambda)\, d_{q} t \\
  \nonumber &-& K(x, x) D_{q^{-1}, t} \phi_{0}(t, \lambda)|_{t=x}+K(x, 0) D_{q^{-1}, t} \phi_{0}(t, \lambda)|_{t=0} \\
   &-& D_{q^{-1}, t} K(x, t)|_{t=0}\phi_{0}(0, \lambda)+D_{q^{-1}, t} K(x, t)|_{t=x}\phi_{0}(x, \lambda).
\end{eqnarray}
Using \eqref{ppp1}, \eqref{l8}, and \eqref{l99} in \eqref{t5}, we get
\begin{eqnarray}\label{tt5}
   \nonumber \lambda \phi(x, \lambda) &=& \lambda \phi_{0}(x, \lambda)+D_{q^{-1}, t} K(x, t)|_{t=x} \phi_{0}(x, \lambda) \\
   &+& \int_{0}^{x}\left(v_{0}(t)K(x, t)-D_{q, t} D_{q^{-1}, t} K(x, t)\right)\phi_{0}(t, \lambda)\, d_{q} t.
\end{eqnarray}
Operating on~\eqref{t7} by $D_{q, x} D_{q^{-1}, x}$, then using \eqref{ff1} and \eqref{ff2}, we obtain
\begin{eqnarray}\label{t6}
\nonumber D_{q, x} D_{q^{-1}, x} \phi(x, \lambda)&=& D_{q, x}D_{{q}^{-1}, x} \phi_{0}(x,\lambda)+
\int_{0}^{x} D_{q, x} D_{q^{-1}, x} K(x, t) \phi_{0}(t, \lambda)\, d_{q} t \\
&+& D_{q, x} \left(K\left(\frac{x}{q}, \frac{x}{q}\right) \phi_{0}(\frac{x}{q},\lambda)\right)+ D_{q, x}  K( x, t)\mid_{t=x} \phi_{0}(x,\lambda).
\end{eqnarray}
Substituting \eqref{p1} and \eqref{l99} into \eqref{t6} yields
 \begin{eqnarray}\label{tt6}
 \nonumber D_{q, x} D_{q^{-1}, x} \phi(x, \lambda)&=& \left(-\lambda +v_{0}(x) +D_{q,  x}  K(x, t)\mid_{t=x} \right) \phi_{0}(x,\lambda)\\
 &+&\int_{0}^{x} D_{q, x} D_{q^{-1}, x} K(x, t) \phi_{0}(t, \lambda)\, d_{q} t.
 \end{eqnarray}
From \eqref{tt5} and \eqref{tt6}, we obtain
\begin{eqnarray}\label{ttt6}
 \nonumber  D_{q, x} D_{q^{-1}, x} \phi(x, \lambda)+\lambda \phi(x, \lambda)&=&\left[ v_{0}(x)+D_{q,x}  K( x, t)\mid_{t=x}+D_{q^{-1}, t} K(x, t)|_{t=x}\right] \phi_{0}(x, \lambda)\\
  \nonumber &+&\int_{0}^{x} \left[ D_{q, x} D_{q^{-1}, x} K(x, t)- D_{q, t} D_{q^{-1}, t} K(x, t)\right] \phi_{0}(t, \lambda)\, d_{q} t\\
   &+& \int_{0}^{x}v_{0}(t) K(x, t)\phi_{0}(t, \lambda)\, d_{q} t.
\end{eqnarray}
Substituting \eqref{l7} and \eqref{t7} into \eqref{ttt6} yields that $\varphi(x, \lambda)$ is a solution of \eqref{slp3}.

Now, we prove \eqref{t78}. From \eqref{t7}, we have
\begin{equation}\label{t877}
  \phi(0,\lambda)=\phi_{0}(0, \lambda)= 1.
\end{equation}
Using \eqref{l99} in \eqref{t7}, we obtain
\begin{align}\label{tt887}
   \nonumber \phi(x, \lambda)=&\phi_{0}(x, \lambda)+\int_{0}^{qx} K(x, t) \phi_{0}(t, \lambda)\, d_{q} t\\
    =&\phi_{0}(x, \lambda)+ q\int_{0}^{x} K(x, qt) \phi_{0}(qt, \lambda)\, d_{q} t.
  \end{align}
  Using \eqref{ff1} in \eqref{tt887} yields
  \begin{equation*}
    \nonumber D_{q^{-1}, x} \phi(x, \lambda)=D_{q^{-1}} \phi_{0}(x, \lambda)+q \int_{0}^{x} D_{q^{-1},x} K(x,qt) \phi_{0}(qt, \lambda) \, d_{q} t+q K(x/q,x) \phi_{0}(x , \lambda).
  \end{equation*}
  Hence,
  \begin{align}\label{t87}
  \nonumber  D_{q^{-1}} \phi(0, \lambda)= & \lim_{n\rightarrow\infty}\left( D_{q^{-1},x} \phi(x, \lambda)\right)(xq^n)\\
  \nonumber =&D_{q^{-1}} \phi_{0}(0, \lambda)+q \phi_{0}(0 , \lambda)\lim_{n\rightarrow\infty}K(xq^{n-1},xq^n) \\
  =& D_{q^{-1}} \phi_{0}(0, \lambda)-q \phi_{0}(0 , \lambda)\lim_{n\rightarrow\infty}F(xq^{n-1},xq^n).
\end{align}
Using \eqref{ppp1} and \eqref{i1} in \eqref{t87}, we get
\begin{equation}\label{t66}
   D_{q^{-1}} \phi(0, \lambda)= h_{0}-q \sum_{n=0}^{\infty} c _{n}=h.
\end{equation}
From \eqref{t877} and \eqref{t66} yields \eqref{t78}. 
\end{proof}
In the following theorems, we assume that the problems $L_q(v_{0},h_{0},H_{0})$ and $L_q(v,h,H)$ have the same spectrum $\{\lambda_n\}^{\infty}_{n=0}.$  
\begin{theorem}
  Let $\{\lambda_{n}\}^{\infty}_{n=0}$ be the eigenvalues of the unperturbed problem~\eqref{p1}-\eqref{ppp1}. Then the eigenfunctions $\{\phi\left(x, \lambda_{n}\right)\}_{n=0}^{\infty}$ defined by equation~\eqref{t7} with $\lambda_{n}$ in place of $\lambda$, can be expressed by the formula
  \begin{equation}\label{t55}
    \phi\left(x, \lambda_{n}\right)=\phi_{0}\left(x, \lambda_{n}\right)-\sum_{k=0}^{\infty} c_{k} \phi\left(x, \lambda_{k}\right) \int_{0}^{x} \phi_{0}\left(t, \lambda_{n}\right) \phi_{0}\left(t, \lambda_{k}\right) d_{q} t,\ x\in A^{*}_{q,a}.
  \end{equation}
  Moreover, the functions $ \{\phi\left(\cdot, \lambda_{n}\right)\}_{n=0}^\infty$ satisfy the boundary condition
\begin{equation*}\label{}
  D_{q^{-1}}\phi\left(a, \lambda_{n}\right)+H \phi\left(a, \lambda_{n}\right)=0,
\end{equation*}
for some real constant $H$.
\end{theorem}
\begin{proof}
 From  \eqref{i2} and \eqref{i1}, we have

\begin{align}\label{t555}
\nonumber K(x, t)=&-F(x, t)-\int_{0}^{x} K(x, s)F(s, t) d_{q}s \\
\nonumber&=-\sum_{k=0}^{\infty} c_{k} \phi_{0}\left(t, \lambda_{k}\right)\left[\phi_{0}\left(x, \lambda_{k}\right)+\int_{0}^{x} K(x, s) \phi_{0}\left(s, \lambda_{k}\right) d_{q}s\right] \\
&=-\sum_{k=0}^{\infty} c_{k} \phi_{0}\left(t, \lambda_{k}\right) \phi\left(x, \lambda_{k}\right).
\end{align}

Substituting \eqref{t555} into \eqref{t7} yields \eqref{t55}. Set $x=a$ in~\eqref{t55}. Then

\begin{align}\label{iso1}
\nonumber \phi\left(a, \lambda_{n}\right)=&\phi_{0}\left(a, \lambda_{n}\right)-\sum_{k=0}^{\infty} c_{k} \phi\left(a, \lambda_{k}\right) \int_{0}^{a} \phi_{0}\left(t, \lambda_{k}\right) \phi_{0}\left(t, \lambda_{n}\right) d_{q} t \\
=&\phi_{0}\left(a, \lambda_{n}\right)-c_{n} \alpha_{n, 0} \phi\left(a, \lambda_{n}\right),
\end{align}
where $\alpha_{n, 0}$ is defined in \eqref{i3}. From \eqref{iso1}, it follows that
\begin{equation}\label{iso2}
  \phi\left(a, \lambda_{n}\right)=\frac{\phi_{0}\left(a, \lambda_{n}\right)}{1+c_{n} \alpha_{n, 0}}.
\end{equation}
The $q$-differentiation of~\eqref{t55} with respect to $D_{q^{-1}, x}$ yields
\begin{align}\label{iso11}
\nonumber D_{q^{-1}, x} \phi \left(x, \lambda_{n}\right)=& D_{q^{-1}, x} \phi_{0}\left(x, \lambda_{n}\right)-\sum_{k=0}^{\infty} c_{k} D_{q^{-1}, x} \phi\left(x, \lambda_{k}\right) \int_{0}^{x} \phi_{0}\left(t, \lambda_{k}\right) \phi_{0}\left(t, \lambda_{n}\right) d_{q} t \\
&-q \sum_{k=0}^{\infty} c_{k} \phi\left(x, \lambda_{k}\right) \phi_{0}\left(x, \lambda_{k}\right) \phi_{0}\left(x, \lambda_{n}\right).
\end{align}
From $\eqref{iso2}$ and the second boundary condition~\eqref{pp1} in \eqref{iso11}, we obtain
\begin{gather}\label{th33}
  \nonumber D_{q^{-1}} \phi \left(a, \lambda_{n}\right)\left(1+c_{n} \alpha_{n, 0}\right)=-H_{0} \phi_{0}\left(a, \lambda_{n}\right)-q \phi_{0}\left(a, \lambda_{n}\right) \sum_{k=0}^{\infty} c_{k} \phi\left(a, \lambda_{k}\right) \phi_{0}\left(a, \lambda_{k}\right) \\
=-H_{0} \phi\left(a, \lambda_{n}\right)\left(1+c_{n} \alpha_{n, 0}\right)-q \left(1+c_{n} \alpha_{n, 0}\right) \phi\left(a, \lambda_{n}\right) \sum_{k=0}^{\infty} c_{k} \frac{\phi_{0}^{2}\left(a, \lambda_{k}\right)}{1+c_{k} \alpha_{k, 0}}.
\end{gather}
Therefore,
\begin{equation*}\label{func1}
  D_{q^{-1}} \phi \left(a, \lambda_{n}\right)=-H \phi \left(a, \lambda_{n}\right),
\end{equation*}
where
\begin{equation}\label{iso3}
  H=H_{0}+ q \sum_{k=0}^{\infty} c_{k} \frac{\phi_{0}^{2}\left(a, \lambda_{k}\right)}{1+c_{k} \alpha_{k, 0}}.
\end{equation}
\end{proof}
\section{A $q$-analog of the Ashrafyan uniqueness theorem}
  In this section, we introduce a $q$-analog of the Ashrafyan uniqueness theorem, namely Theorem C. We present some definitions that will guide our next theorems, aiding in the proof of the Ashrafyan uniqueness theorem.

 Let $f, g$ be entire functions and $a \in \mathbb{C}$, we say that \begin{equation*}
   f(z)=O(g(z)), \quad \text { as } \quad z \rightarrow a,
 \end{equation*}
if $f(z) / g(z)$ is bounded in a neighborhood of $a$. 
\begin{equation*}
  f(z) \sim g(z), \text { as } z \rightarrow a, \text { if } \lim _{z \rightarrow a} \frac{f(z)}{g(z)}=1.
\end{equation*}
If $f(z):=\sum_{n=0}^{\infty} a_n z^n$ is an entire function, then the maximum and  minimum modulus are defined respectively for $r>0$ by
\begin{equation*}\label{}
  M(r ; f):=\sup \{|f(z)|:|z|=r\},
\end{equation*}
\begin{equation*}\label{}
  \mu(r ; f):=\inf \{|f(z)|:|z|=r\}.
\end{equation*}
The order  $\rho(f)$ of  $f$ is defined as 
\begin{equation*}
  \rho(f):=\displaystyle {\limsup _{r \rightarrow \infty}} \frac{\log \log M(r, f)}{\log r}=\limsup _{n \rightarrow \infty} \frac{n \log n}{-\log |a_n|};
\end{equation*}
see~\cite{boas11}. In the following, the functions $\phi_1(\cdot, \lambda)$ and $\phi_2(\cdot, \lambda)$ are the solutions of \eqref{slp3} with the initial conditions
 \begin{equation*}\label{z1}
   D_{q^{-1}}^{(j-1)} \phi_i(0, \lambda)=\delta_{i j}, \quad 1 \leqslant i, j \leqslant 2, \quad \lambda \in \mathbb{C},
 \end{equation*}
where $D_q^{(0)}$ is the identity operator. In this case, any solution of \eqref{slp3} is a linear combination of $\left\{\phi_{1}, \phi_{2}\right\}$.   In \cite{annaby11}, the authors defined the characteristic determinant $\Delta(\lambda)$ as 
\begin{equation}\label{cf2}
  \Delta(\lambda)=U_1\left(\phi_1\right) U_2\left(\phi_2\right)-U_1\left(\phi_2\right) U_2\left(\phi_1\right),
\end{equation}
where $U_1$ and $U_2$ are defined  in \eqref{slp4} and \eqref{slp5}, respectively. 

 In their paper \cite{Annaby07}, Annaby and Mansour presented the asymptotic behavior of the functions  $\cos (z ; q)$ and $\sin (z ; q)$ by employing the technique developed by Bergweiler and Hayman to investigate the asymptotic properties of the solutions of functional equations, see~\cite{BH,Bergweiler19}. 
 
 \begin{theorem}\label{asy333} \textup{\cite{Annaby07}}
   As $n\rightarrow \infty$, we have
\begin{align*}
 & x_n=\frac{q^{-n+\frac{1}{ 2}}}{(1-q)}\left(1+O\left(q^n\right)\right), \\
& y_n=\frac{q^{-n}}{(1-q)}\left(1+O\left(q^n\right)\right),
\end{align*}
where $\{x_n \}$ and $\{ y_n \}$ are the  positive zeros of $\cos (\cdot ; q)$ and $\sin (\cdot ; q)$, respectively.
 \end{theorem}
 
 The authors in \cite{annaby11} examined the asymptotic formulas for eigenvalues and eigenfunctions of the q-Sturm--Liouville problem~\eqref{slp3}-\eqref{slp5} and established  the following theorems:
 
\begin{theorem}\label{asy2} \textup{\cite{annaby11}}
   Let $\lambda$ be a complex number, and let $s=\sqrt{\lambda}$. Then the asymptotic formulaes
\begin{equation}\label{as1}
  \phi_1(x, \lambda)=\cos (s x ; q)+O\left(|s|^{-1}exp\left(\dfrac{-\left(\log |s| x (1-q)\right)^2}{\log q}\right)\right),
\end{equation}
\begin{equation}\label{as2}
\phi_2(x, \lambda)=\frac{\sin (s x ; q)}{s}+O\left(|s|^{-2}exp\left(\dfrac{-\left(\log |s| x (1-q)\right)^2}{\log q}\right)\right),
\end{equation}
 holds as $|\lambda| \longrightarrow \infty$, where for each $x \in (0,a]$ the $O$-terms are uniform on $\{xq^n,\ n\in \mathbb{N}_0\}$. Moreover, if $v(\cdot)$ is bounded on $[0,a]$, the $O$-terms~\eqref{as1}-\eqref{as2} will be uniform for all $x\in [0, a]$.
\end{theorem}

\begin{theorem}\label{asy22}\textup{\cite{annaby11}}
  As $|\lambda| \longrightarrow \infty$,  the asymptotic relation of $\Delta(\lambda)$ defined in \eqref{cf2} has the following form:
\begin{equation}\label{asy33}
  \Delta(\lambda) =-a_{12} a_{22}\sqrt{q} s \sin(s q^{-1 / 2} a ; q)+O\left(\exp \left(\frac{-\left(\log |s| q^{-1 / 2} a(1-q)\right)^2}{\log q}\right)\right),\quad a_{12} a_{22}\neq0.
\end{equation}
\end{theorem}

\begin{theorem}\label{thas}\textup{\cite{annaby11}}
  The positive zeros $\left\{\lambda_n\right\}$ of $\Delta(\lambda)$ that defined in \eqref{cf2} are given asymptotically as $n \longrightarrow \infty$ by
\begin{equation*}\label{}
  \lambda_n= \begin{cases}\frac{q}{a^2} y_n^2\left(1+O\left(q^{n / 2}\right)\right), & a_{12} a_{22} \neq 0, \\ \frac{q}{a^2} x_n^2\left(1+O\left(q^{n/ 2}\right)\right), & a_{12}=0, \\ \frac{1}{a^2} x_n^2\left(1+O\left(q^{n / 2}\right)\right), & a_{22}=0,\end{cases}
\end{equation*}
where $\{x_n\}, \{y_n \}$ represent  the positive zeros of $\cos(\cdot ; q)$ and $\sin(\cdot ; q)$,  respectively.
\end{theorem}

From Theorem~\ref{asy333} and Theorem~\ref{thas}, we obtain the following corollary.
\begin{corollary}\label{con9}\textup{\cite{annaby11}}
  As $n \longrightarrow \infty$,
 \begin{equation}\label{asy3}
  s_n:= \sqrt{\lambda_n}= \begin{cases}\frac{q^{-n+1 / 2}}{a(1-q)}\left(1+O\left(q^{n / 2}\right)\right), & a_{12} \neq 0, \\ \frac{q^{-n+1}}{a(1-q)}\left(1+O\left(q^{n / 2}\right)\right), & a_{12}=0 .\end{cases}
 \end{equation}
\end{corollary}

\begin{prop}\label{prop2}
  If $\{y(\cdot,\lambda_n)\}_{n=0}^\infty$ is the set of solutions of the problem~\eqref{slp3}-\eqref{slp5}, then there exist positive constants $k_{1}$, $ k_{2}$ (independent of $n$) such that
  \begin{equation}\label{prop22}
    k_{1} \exp\left({\dfrac{-(\log^2 |s_n| a(1-q))}{\log q}}\right) \leq |y(a,\lambda_n)| \leq k_{2} \exp\left({\dfrac{-(\log^2 |s_n| a(1-q))}{\log q}}\right),\quad n\in \N.
  \end{equation} 
\end{prop}
\begin{proof}
  Since $\{y(x,\lambda_n)\}$ is a solution of the problem~\eqref{slp3}-\eqref{slp5}, then $y(x,\lambda_n)$ satisfies the functional equation~\eqref{fff3}. Furthermore, $y(a,\lambda)$ is an entire function in $\lambda $ of order zero; see \cite{boas11}. From Remark~\ref{eu1}, we have $y(a,\lambda_n)\neq 0$ for all $n\in \N$.  Let   the set $E$  be defined as 
  \begin{equation*}
    E:=\R \setminus \{ t_n\in \R : y(a,t_n)=0 \}.
  \end{equation*}
  Hence,
   \begin{equation*}
     \lim _{n \rightarrow \infty} \dfrac{m\{E \cap(0, \lambda_n)\}}{\lambda_n}=1,
   \end{equation*}
  From \cite[p.1]{kubota69}, we have 
   \begin{equation*}
    M\left( \lambda_n; y(a,\lambda)\right) \sim   \mu\left( \lambda_n; y(a,\lambda)\right)\; \text{ as}\  n\rightarrow \infty.
   \end{equation*} 
   Since \[|y(a,\lambda)|\sim exp \left(\dfrac{-(\log^2 |s| a(1-q))}{\log q}\right)\; \mbox{for sufficiently large}\; \lambda,\]
   where $s:=\sqrt{\lambda}$,  see\cite{annaby11}. Hence, there exist positive constants $k_{1}, k_{2}$ such that \eqref{prop22} holds. 
\end{proof}

\begin{corollary}\label{phic}
 Let  $\{\phi\left(\cdot, \lambda_n \right)\}$   and $\{\phi_0\left(\cdot, \lambda_{n,0} \right)\}$ be  the eigenfunctions of the $q$-Sturm--Liouville problems   \eqref{slp3}--\eqref{is4} and  \eqref{ppp1}--\eqref{pppp1},  respectively. Then  we have the asymptotic formulae as $n\rightarrow\infty$
 \begin{equation}\label{phi1}
  \phi(x, \lambda_n)=\cos (s_n x ; q)+\frac{h\sin(s_n x ; q)}{s_n}+O\left(|s_n|^{-1}exp\left(\dfrac{-\left(\log |s_n| x (1-q)\right)^2}{\log q}\right)\right),\ 
\end{equation}
  \begin{equation}\label{phi2}
  \phi_0(x, \lambda_{n,0})=\cos (s_{n,0} x ; q)+\frac{h_0 \sin(s_{n,0} x ; q)}{s_{n,0}}+O\left(|s_{n,0}|^{-1}exp\left(\dfrac{-\left(\log |s_{n,0}| x (1-q)\right)^2}{\log q}\right)\right),
\end{equation}
 where $s_n=\sqrt{\lambda_n},\ s_{n,0}=\sqrt{\lambda_{n,0}}$ and for each $x \in (0,a]$ the $O$-terms are uniform on $\{xq^n,\ n\in \mathbb{N}_0\}$. Moreover, if $v(\cdot)$ is bounded on $[0,a]$, the $O$-terms~\eqref{phi1}-\eqref{phi2} will be uniform for all $x\in [0, a]$..  
\end{corollary}

\begin{proof}
  Since $\phi\left(x, \lambda\right)$ is the solution of $\eqref{slp3}$ with the boundary condition~\eqref{is4}, then
\begin{equation}\label{d55}
  \phi\left(x, \lambda\right)=\phi_{1}\left(x, \lambda\right)+h \phi_{2}\left(x, \lambda\right),\quad x \in (0,a].
\end{equation}
Using Theorem~\ref{asy2}, we obtain \eqref{phi1}. In the same way, we prove \eqref{phi2}.
\end{proof}

\begin{remark}\label{ordz}
From Corollary~\ref{phic}, $M(r ; \phi) \sim exp \left(\dfrac{-\left(\log r(1-q)\right)^2}{\log q}\right)$ as $r\rightarrow \infty.$
Hence, $\phi$ is an entire function of order zero, and so $\Delta(\lambda)$; see \textup{\cite{boas11}}.  
\end{remark}

\begin{lemma}\label{lem1}
 If $\{ \phi\left(\cdot, \lambda_n\right) \}$ denotes the eigenfunctions corresponding to the spectrum $\{\lambda_n\}$ of the problem $L_q(v,h,H)$, then the characteristic determinant $\Delta(\lambda)$ and the norming constants $\{\alpha_{n}\}$ can be represented as follows:
 \begin{eqnarray}\label{dd}
   \Delta(\lambda)  &=& -D_{q^{-1}} \phi(a, \lambda)-H \varphi(a, \lambda),
 \end{eqnarray}
\begin{eqnarray*}
   \alpha_{n} &=& \phi\left(a, \lambda_{n}\right) \dot{\Delta}\left(\lambda_{n}\right), \quad n\in \N_{0},
\end{eqnarray*}
where the dot of $\Delta(\lambda)$ is the derivative with respect to $\lambda$.
 \end{lemma}
 
 \begin{proof}
Substituting \eqref{slp4} and \eqref{slp5} into \eqref{cf2}, we have 
\begin{equation*}\label{}
   \Delta(\lambda)= a_{11}\left[ a_{21}\phi_{2}\left(a, \lambda\right)+a_{22} D_{q^{-1}} \phi_{2}\left(a, \lambda\right)\right]- a_{12}\left[a_{21}\phi_{1}\left(a, \lambda\right)+a_{22} D_{q^{-1}} \phi_{1}\left(a, \lambda\right)\right].
 \end{equation*}
 Comparing \eqref{is4} with  \eqref{slp4} and  \eqref{slp5}, we obtain
 \[a_{11}=-h,\ a_{12}=a_{22}=1,\ a_{21}=H.\] Consequently,
 \begin{equation}\label{d44}
   \Delta(\lambda)=-H\left[\phi_{1}\left(a, \lambda\right)+h \phi_{2}\left(a, \lambda\right)\right]- \left[D_{q^{-1}} \phi_{1}\left(a, \lambda\right)+h\phi_{1}\left(a, \lambda\right)\right].
 \end{equation}
Using \eqref{d55} in \eqref{d44} yields \eqref{dd}.
\par  For arbitrary $\lambda$ and $\mu$,  set $y=\phi(x, \lambda)$ and $ z=\phi(x, \mu)$. Hence,   from  the Green's identity \eqref{lin}, we obtain
 \begin{equation*}\label{eqr1}
   (\mu-\lambda) \int_{0}^{a} y(x)z(x)\, d_{q} x= W_{q^{-1}}\left(y, z\right)(0)- W_{q^{-1}}\left(y, z\right)(a).
 \end{equation*}
 Using \eqref{wr} and the initial condition in \eqref{is4} yields 
 \begin{equation*}\label{}
   \int_{0}^{a} y(x)z(x) \,d_{q} x =\dfrac{\phi(a, \mu) D_{q^{-1}} \phi(a, \lambda)-\phi(a, \lambda) D_{q^{-1}} \phi(a, \mu)}{(\mu-\lambda)}.
 \end{equation*}
 Taking  the limit as $\mu$ approaches  $\lambda $, yields 
 \begin{equation*}\label{}
    \int_{0}^{a} \phi^{2}(x, \lambda) d_{q} x = \dot{\phi}(a, \lambda) D_{q^{-1}} \phi(a, \lambda )-\phi(a, \lambda ) D_{q^{-1}} \dot{\phi}(a, \lambda ).
 \end{equation*}
 Setting $\lambda=\lambda_n$ and using the boundary condition in \eqref{is4}, we obtain
  \begin{equation*}\label{}
  \alpha_{n}=\phi\left(a, \lambda_{n}\right) \left[-H \dot{\phi}(a, \lambda_n)-D_{q^{-1}} \dot{\phi}(a,\lambda_n)\right]= \phi\left(a, \lambda_{n}\right) 
   \dot{\Delta}\left(\lambda_{n}\right), \quad n\in \N_{0}.
  \end{equation*}
\end{proof}

\begin{theorem}\label{eql}
  The specification of the spectrum $\left\{\lambda_n\right\}^\infty_{n=0}$ of the problem~$L_q(v,h,H)$ uniquely determines the characteristic function $\Delta(\lambda)$ by the formula
  \begin{equation}\label{gen1}
    \Delta(\lambda)=a \left(\lambda_0-\lambda\right) \prod_{n=1}^{\infty} \frac{\lambda_n-\lambda}{\frac{q y_n^2}{a^2}}.
  \end{equation}
\end{theorem}

\begin{proof}
It follows from Remark~\ref{ordz} that $\Delta(\lambda)$ is an entire  function in $\lambda$ of order zero.  Consequently by Hadamard's factorization theorem \cite[chapter~5]{Stein03}, $\Delta(\lambda)$ has the following form
  \begin{equation}\label{gen}
    \Delta(\lambda)=C \prod_{n=0}^{\infty}\left(1-\frac{\lambda}{\lambda_n}\right),
  \end{equation}
  where $C$ is a constant. Consider the function
  \begin{equation}\label{}
    \tilde{\Delta}(\lambda):=-s \sqrt{q} \sin( s a q^{-\frac{1}{2}};q)=-\lambda a \prod_{n=1}^{\infty}\left(1-\dfrac{a^2\lambda}{q y_n^2}\right),
  \end{equation}
  where $\{y_n\}$ are the positive zeroes of $\sin(\cdot;q)$. Then
  \begin{equation*}\label{}
    \dfrac{\Delta(\lambda)}{\tilde{\Delta}(\lambda)}=C \dfrac{\lambda-\lambda_0}{a \lambda_0 \lambda}\;\;\; \prod_{n=1}^{\infty} \dfrac{qy^2}{a^2\lambda_n} \prod_{n=1}^{\infty}\left(1+\dfrac{\lambda_n-\frac{q y_n^2}{a^2}}{\frac{q y_n^2}{a^2}-\lambda}\right).
  \end{equation*}
From Corollary~\ref{con9} and \eqref{asy33}, we conclude
\begin{equation*}
  \lim _{\lambda \rightarrow-\infty} \frac{\Delta(\lambda)}{\tilde{\Delta}(\lambda)}=1, \quad \lim _{\lambda \rightarrow-\infty} \prod_{n=1}^{\infty}\left(1+\dfrac{\lambda_n-\frac{q y_n^2}{a^2}}{\frac{q y_n^2}{a^2}-\lambda}\right)=1.
\end{equation*}
Hence,
\begin{equation}\label{const}
  C=a \lambda_0 \prod_{n=1}^{\infty} \dfrac{a^2\lambda_n}{q y_n^2}.
\end{equation}
Substituting \eqref{const} into \eqref{gen} yields \eqref{gen1}, which completes the proof.

It is worth noting that from Theorem~\ref{thas}, 
\begin{equation*}
  \lambda_n\simeq \dfrac{q}{a^2} y_n^2\left(1+O\left(q^{n / 2}\right)\right),\ \text{as}\  n\rightarrow \infty.
\end{equation*}
Hence,
\begin{equation*}
\prod_{n=1}^{\infty} \dfrac{\lambda_n-\lambda}{q y_n^2/a^2}=\prod_{n=1}^{\infty}\left(1+r_n-\dfrac{a^2 \lambda}{q y_n}\right),\  r_n=O\left(q^{n / 2}\right) \text{as}\  n\rightarrow \infty.
 \end{equation*}
Therefore, the infinite product in \eqref{gen1} is convergent. Similarly, we can prove that the infinite product in \eqref{const} is convergent.
\end{proof}
\begin{prop}\label{alpha} Let $\{\lambda_{n}, \alpha_{n} \}_{n=0}^\infty$ and $\{\lambda_{n,0}, \alpha_{n,0}\}_{n=0}^\infty$ be the spectral data of the problems
 $L_q(v,h,H)$ and $L_q(v_0,h_0,H_0)$, respectively.
 If $\lambda_n= \lambda_{n,0}$, then
\begin{equation}\label{inv22}
    \dfrac{\alpha_{n}}{\alpha_{n,0}}=\dfrac{\phi(a, \lambda_n)}{\phi_0(a, \lambda_n)}.
\end{equation}
\end{prop}
\begin{proof}
  Using Lemma~\ref{lem1}, we have
  \begin{equation*}\label{alpha2}
    \dfrac{\alpha_{n}}{\alpha_{n,0}}=\dfrac{\phi(a, \lambda_n)\dot{\Delta}\left(\lambda_{n}\right)}{\phi_0(a, \lambda_{n,0})\dot{\Delta_0}\left(\lambda_{n,0}\right)},
\end{equation*}
where $\Delta\left(\lambda_{n}\right),$ $\Delta_0\left(\lambda_{n,0}\right)$ are the characteristic functions associated with the problems $L_q(v,h,H))$ and $L_q(v_0,h_0,H_0)$, respectively. 
Since $\lambda_n= \lambda_{n,0}$ for all $n\in\N_0$, then by using Theorem~\ref{eql}, we have $\Delta\left(\lambda_{n}\right)=\Delta_0\left(\lambda_{n,0}\right)$ for all $\lambda \in\C.$ Consequently, \[\dot{\Delta}\left(\lambda_{n}\right)=\dot{\Delta}_0\left(\lambda_{n,0}\right)\; \; (n\in\mathbb{N}_0), \]  and \eqref{inv22} follows.
\end{proof}

\begin{theorem}
Let $\{\phi_0\left(x, \lambda_{n,0} \right)\}_{n=0}^{\infty}$ and $\{\phi\left(x, \lambda_{n} \right)\}_{n=0}^{\infty}$ be the eigenfunctions of the problems $L_q(v_0,h_0,H_0),\ L_q(v,h,H)$, respectively. If 
\[  h=h_0,\;       \lambda_{n}=\lambda_{n,0}\qquad \left(n \in \N_{0}\right),\]  then the series 
\begin{equation}\label{ill1}
F(x, y)=\sum_{n=0}^{\infty} (\dfrac{1}{\alpha_{n}}-\frac{1}{\alpha_{n, 0}})\; \phi_{0}\left(x, \lambda_{n,0}\right) \phi_{0}\left(y, \lambda_{n,0}\right),
\end{equation}
converges uniformly for $(x,t)\in A_{q,a}\times A_{q,a}$.
\end{theorem}
\begin{proof}
Since $\lambda_{n}=\lambda_{n,0}$, then equation \eqref{ill1} can be written as 
\begin{equation}\label{cov8}
F(x, y)=\sum_{n=0}^{\infty} \dfrac{1}{\alpha_{n}} (1-\dfrac{\alpha_n}{\alpha_{n, 0}}) \phi_{0} \left(x, \lambda_{n}\right) \phi_{0} \left(y, \lambda_{n}\right).
\end{equation}
Using~\eqref{inv22} and Corollary~\ref{phic} at $h=h_0,$ we conclude that
\begin{equation}\label{conv1}
  \dfrac{\alpha_n}{\alpha_{n, 0}}= \dfrac{\phi \left(a, \lambda_{n}\right)}{\phi_{0} \left(a, \lambda_{n}\right)} =1+O(\frac{1}{s_n}) \ \text{as}\  n \rightarrow \infty.
\end{equation}
From \eqref{conv1}, there exists  a positive constant $k$ such that 
\begin{equation}\label{conv2}
 |1- \dfrac{\alpha_n}{\alpha_{n, 0}}| \leq  \dfrac{k}{s_n},\quad  (n \in \N_{0}).
\end{equation}
From \eqref{is333} and \eqref{int}, we have
\begin{equation}\label{conv4}
  \alpha_n \geq a(1-q) \phi^2 \left(a, \lambda_{n}\right).
\end{equation}
Applying Proposition~\ref{prop2}, there exist a positive constant $k_1$ such that
 \begin{equation}\label{conv5}
   \dfrac{1}{|\alpha_n |}\leq \dfrac{1}{a(1-q){k_1}^2 \exp\left({\dfrac{-2(\log^2 |s_n| a(1-q))}{\log q}}\right)}\quad ( n\in \N_{0}).
 \end{equation}
 Also, from Remark~\ref{ordz},
\begin{equation}\label{conv6}
  |\phi_{0} \left(x, \lambda_{n}\right) \phi_{0} \left(y, \lambda_{n}\right)| \leq {k_2}^2 \exp\left({\dfrac{-2(\log^2 |s_n| a(1-q))}{\log q}}\right)\quad ( n\in \N_{0}),
\end{equation}
where $k_2$ is a positive constant.
From \eqref{conv5}, \eqref{conv6}, we have
\begin{equation}\label{conv7}
  \dfrac{|\phi_{0} \left(x, \lambda_{n}\right) \phi_{0} \left(y, \lambda_{n}\right)|}{|\alpha_n|} \leq M\quad ( n\in \N_{0}),
\end{equation}
where $M=\dfrac{ k^2_2 }{a(1-q){k_1}^2}.$
Substituting~\eqref{conv2} and  \eqref{conv7} into \eqref{cov8}, we obtain
\begin{equation}\label{cov8}
|F(x, y)| = \sum_{n=0}^{\infty} \dfrac{1}{|\alpha_{n}|} |(1-\dfrac{\alpha_n}{\alpha_{n, 0}})| | \phi_{0} \left(x, \lambda_{n}\right) \phi_{0} \left(y, \lambda_{n}\right) | \leq M k\sum_{n=0}^{\infty} \dfrac{1}{|s_n|}.
\end{equation}
From  the asymptotic of the sequence $\{s_n\}$, see~\eqref{asy3}, the series $\sum_{n=0}^{\infty} \dfrac{1}{|s_n|}$ is convergent. This is the final step of the proof.
\end{proof}

 Now, we introduce a $q$-analog of Theorem C, the  Levinson-Marchenko theorem.

 \begin{theorem}\label{nt80}
  Let $\{\lambda_{n}, \alpha_{n} \}_{n=0}^{\infty}$ and $\{{\lambda}_{n,0}, \alpha_{n,0}\}_{n=0}^{\infty}$ be the spectral data of the problems
 $L_q(v,h,H)$ and $ L_q(v_{0},h_{0},H_{0}),$ respectively. If $h=h_0,$ ${\lambda}_{n}={\lambda}_{n,0}$ and 
 $\alpha_{n}\geq\alpha_{n,0}$, $n \in \N_{0}$, then $L_q=L_{q,0}$, i.e., $v(x)=v_{0}(x)$ on $A_{q,a}$ (qa.e.) and $H=H_{0}$.
\end{theorem}

\begin{proof}
Let us consider problem~\eqref{p1} with the boundary conditions \eqref{pp1} and problem \eqref{slp3} with the boundary conditions \eqref{is4}. These problems have the following norming constants:
\[
  \alpha_{n, 0}=\int_{0}^{a} \phi_{0}^{2}\left(x, \lambda_{n,0}\right) d_{q} x, \quad \alpha_{n}=\int_{0}^{a} \phi^{2}\left(x, \lambda_{n}\right) d_{q} x.
\]
The $q$-Gelfand--Levitan operator~\eqref{i2}  establishes the connection between  these two problems as stated  on Theorem~\ref{nt}.  Substituting~\eqref{inv22} into \eqref{iso2} yields
\begin{equation}\label{nt81}
  c_{n}=\left(\frac{1}{\alpha_{n}}-\frac{1}{\alpha_{n, 0}}\right).
\end{equation}
 Using \eqref{nt81} in \eqref{i1} at ${\lambda}_{n}={\lambda}_{n,0}$, we get
\begin{equation*}
F(x, y)=\sum_{n=0}^{\infty} (\dfrac{1}{\alpha_{n}}-\frac{1}{\alpha_{n, 0}}) \phi_{0}\left(x, \lambda_{n}\right) \phi_{0}\left(y, \lambda_{n}\right).
\end{equation*}
From \eqref{nt81}, \eqref{t78},  and \eqref{iso3}, we obtain
\begin{equation}\label{nt82}
 h= h_{0}-q \sum_{n=0}^{\infty} \left(\frac{1}{\alpha_{n}}-\frac{1}{\alpha_{n, 0}}\right),
\end{equation}
and
\begin{equation*}\label{nt84}
   H=H_{0}+ q \sum_{n=0}^{\infty} \dfrac{\left(\alpha_{n, 0}-\alpha_{n}\right)}{\alpha_{n, 0}}\phi_{0}^{2}\left(a, \lambda_{n}\right).
\end{equation*}
Since $h=h_{0}$, then \eqref{nt82} yields
\begin{equation}\label{nt85}
 \sum_{n=0}^{\infty} \left(\frac{1}{\alpha_{n}}-\frac{1}{\alpha_{n, 0}}\right)=0.
\end{equation}
Since $\alpha_{n} \geq \alpha_{n, 0}$ for all $n \in \N_{0}$, then from equation \eqref{nt85}, we conclude that $\alpha_{n} = \alpha_{n, 0}$ for all $n \in \N_{0}$. Therefore, using Theorem~\ref{nt60}, we obtain $v(x)= v_{0}(x)$ on $A_{q,a}$  (qa.e.) and $H=H_{0}$. This completes the proof.
\end{proof}

\noindent {\bf Availability of data and materials}

The data and material in this paper are original.

\vskip .5 cm
\noindent {\bf Competing interests} 

The authors declare that they have no competing interests.

\vskip .5 cm
\noindent {\bf Funding} 

Not applicable.
\vskip .5 cm

\bibliographystyle{plain}
\bibliography{Ref}
\end{document}